\documentclass[10pt]{iopart}
\usepackage{iopams}
\usepackage{epsfig}

\usepackage{epsfig}

\usepackage{xcolor}
\usepackage[active]{srcltx}
\usepackage{psfrag}

\begin{document}

\newtheorem{lm}{Lemma}
\newtheorem{theorem}{Theorem}
\newtheorem{crl}{Corollary}
\newtheorem{prop}{Proposition}

\newtheorem{df}{Definition}
\newtheorem{remark}{Remark}

\newcommand{\vv}{\mbox{{\bf {\em v}}}}   % bold vector $v$
\newcommand{\vw}{\mbox{{\bf {\em w}}}}   % bold vector $w$
\newcommand{\vr}{\mbox{{\bf {\em r}}}}   % bold vector $r$
\newcommand{\vM}{\mbox{{\bf {\em M}}}}   % bold vector $M$
\newcommand{\vo}{\boldsymbol \omega}   % bold vector $omega$
\newcommand{\vg}{\boldsymbol \gamma}   % bold vector $gamma$
\newcommand{\vI}{\mbox{{\bf {\em I}}}}   % bold vector $I$
\newcommand{\va}{\mbox{{\bf {\em a}}}}   % bold vector $a$
\newcommand{\vb}{\mbox{{\bf {\em b}}}}   % bold vector $b$
\newcommand{\vc}{\mbox{{\bf {\em c}}}}   % bold vector $c$

\eqnobysec

\def\fl{\!}
\def\TT{{\cal T}}
\def\ds{\displaystyle}
\def\proof{\noindent {\sf Proof.}  }
\def\qed{\hfill $\Box$ \\ \bigskip}
\def\Id{\mbox{Id}}
\def\loc{\mbox{loc}}

\newcommand{\N}{\ensuremath{\mathbb{N}}}
\newcommand{\Z}{\ensuremath{\mathbb{Z}}}
\newcommand{\Q}{\ensuremath{\mathbb{Q}}}
\newcommand{\R}{\ensuremath{\mathbb{R}}}
\newcommand{\C}{\ensuremath{\mathbb{C}}}
\newcommand{\SI}{\ensuremath{\mathbb{S}^1}}
\newcommand{\T}{\ensuremath{\mathbb{T}}}

% macros

\newcommand{\Rrev}{R}   % reversor (nonlinear in general)
\newcommand{\Lrev}{L}   % linear reversor (linear part of R)
\newcommand{\Wu}{W^u}  % unstable invariant manifold
\newcommand{\Ws}{W^s}  % stable invariant manifold
\newcommand{\sech}{\mathrm{sech}}

\definecolor{Sepia}{rgb}{0.6,0.3,0.1}
%\definecolor{Sepia}{rgb}{0.6,0.3,0.1}
\definecolor{RawSienna}{rgb}{0.2,0.7,0.1}
\definecolor{NavyBlue}{rgb}{0.2,0,0.5}
\definecolor{PineGreen}{rgb}{0,0.5,0}
\definecolor{Sepia1}{rgb}{0.9,0.0,0.1}
\definecolor{Kirpich}{rgb}{0.8,0.1,0.1}
\definecolor{Sepia2}{rgb}{0.4,0.,0.4}
\definecolor{Sepia3}{rgb}{0.4,0.,0.}

\newcommand{\rk}[1]{\textcolor{red}{#1}}
\newcommand{\yel}[1]{\textcolor{yellow}{#1}}

% Tomas (macros defined on May 13th)
\newcommand{\Mt}{\tilde{M}}
\newcommand{\ct}{\tilde{c}}

\hyphenation{dif-feo-mor-phisms}

% %%%%%%%%

\title[Lorenz-like attractors in  nonholonomic models of Celtic stone]
{Lorenz-like attractors in  nonholonomic models of Celtic stone.}

\author{Gonchenko A.S.\dag and Gonchenko S.V.\dag}

\address{\dag\ Lobachevsky State University of Nizhni Novgorod, Nizhni Novgorod, Russia}

%\address{\ddag\ Institute for Applied Mathematics \& Cybernetics, N.Novgorod, Russia}

%\theoremstyle{plain}

%\begin{document}

\vspace{0.5cm}
%\begin{center}

\begin{abstract}{We study chaotic dynamics in  nonholonomic model of Celtic stone.
We show that, for certain values of parameters characterizing geometrical and physical properties of the stone, a strange Lorenz-like attractor is observed in the model. We study also bifurcation scenarios for appearance and break-down of this attractor.
}
\end{abstract}

\vspace{5mm}

\section{Introduction}\label{sec:intro}

It is well known that the dynamics of three-dimensional systems of differential equations
can be drastically different from that of two-dimensional systems. The same situation takes place for the dynamics of three-dimensional diffeomorphisms in contrast with
the dynamics of two-dimensional ones. Though chaotic behavior can be observed
both in two-dimensional and three-dimensional maps, nevertheless, especially when the Jacobian of a three-dimensional
map is not too close to zero, the chaos can take forms which are totally different from ones observed
in two-dimensional diffeomorphisms.

One of these phenomena, discovered in \cite{GOST05}, is the emergence of Lorenz-like attractors
for diffeomorphisms of dimension 3 and higher. We will call such attractors {\em discrete Lorenz attractors}. The simplest model of such attractor can be given by the attractor
in the Poincar\'e map
%for period perturbation
for periodically perturbed three-dimensional flow having the classical Lorenz attractor. The latter
is a genuine strange attractor \cite{ABS77,ABS83}: it contains no stable periodic orbits,
and every orbit in such attractor has positive maximal Lyapunov exponent. Moreover, these properties are
robust,
%(they persist at small changes of the right-hand sides)
even though the attractor itself is structurally
unstable \cite{GW79,ABS83}.\footnote{We note that this property seemingly does not hold for many ``physical'' attractors observed in numerical experiments, where an observed chaotic behavior can easily correspond to some periodic orbit with a very large period (plus inevitable noise); see more discussion in \cite{N74,AfrSh83}. In particular, H\'enon-like strange attractors, \cite{BC,MV93}, that occur often in two-dimensional maps may transform into stable long-period orbits by arbitrarily small changes of parameters \cite{Ures95}.}
When the perturbation is small the obtained discrete Lorenz attractor possesses the same properties \cite{TS08}.

Note that the dynamics of discrete Lorenz attractors is actually different from the dynamics of Lorenz attractors
for flows (by taking the flow suspension of a three-dimensional diffeomorphism with a Lorenz-like attractor we obtain a four-dimensional
flow whose attractor may contain Newhouse wild sets and saddle periodic orbits with three-dimensional unstable manifolds \cite{TS08},
as well as saddle invariant tori and heterodimensional cycles; none of these objects exists within the classical Lorenz attractor).
However, visually, a Lorenz-like attractor of a diffeomorphism may look quite similar to the classical
Lorenz attractor, see e.g. Fig.~\ref{fig3-4}.

The reason for the robust chaoticity of both the flow and discrete Lorenz attractors is that they possess a {\em pseudo-hyperbolic}
structure, \cite{TS98,TS08}
(see the exact Definition 1 in Section 3.1).

As a result, the discrete Lorenz attractor does not contain stable periodic orbits,
 and moreover, no any stable orbit can emerge in its neighbourhood for small perturbations of the diffeomorphism.
Thus, the fact that these attractors exist in important class of three-dimensional maps (see \cite{GOST05,GGOT13}) is remarkable. Moreover in \cite{GGS12,GGKT14}, it was proposed a simple and universal scenario of emergence of discrete Lorenz attractors in 3D maps; this scenario is the result of the bifurcation chain which includes the sequence
$$
\mbox{periodic attractor}\;\;\Rightarrow\;\;\mbox{period doubling}\;\;\Rightarrow\;\;\mbox{homoclinic attractor},
$$
\begin{figure}[htb]
\centerline{
  \includegraphics[width=14cm]{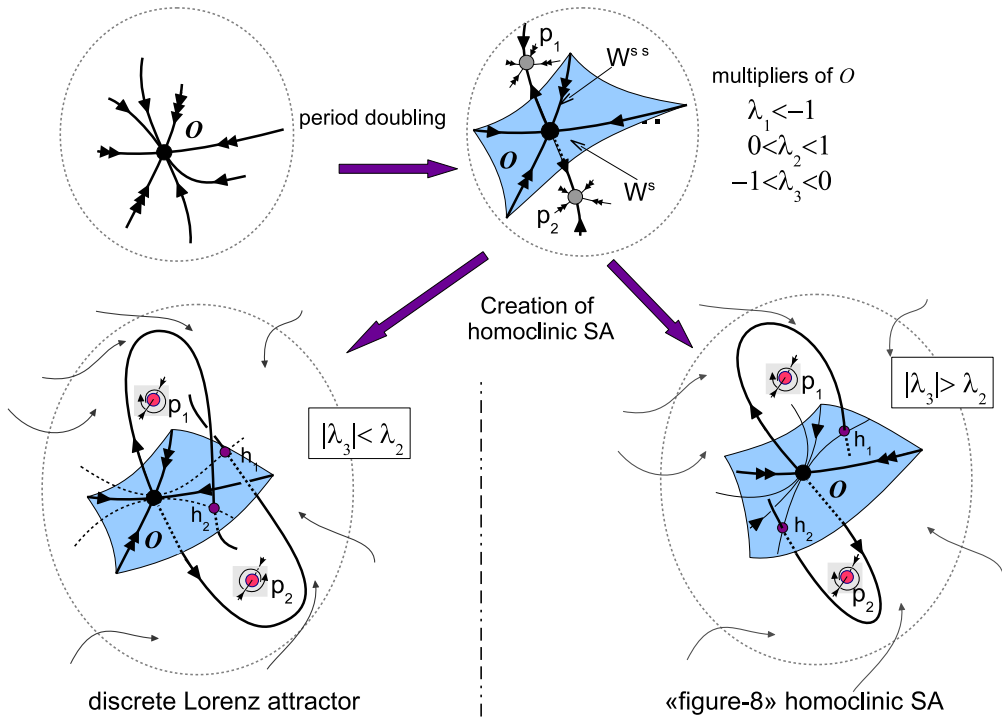}}
  \caption{{\footnotesize Main steps of creation of a homoclinic strange attractor: a discrete Lorenz attractor (below left), a ``figure-eight'' pseudo-hyperbolic strange attractor (below right).}}\label{Lor8attr}
\end{figure}
In Fig.~\ref{Lor8attr} this scenario is illustrated for the case when the periodic attractor $O$ is an asymptotically stable fixed fixed point; then
%$\Rightarrow$
 this point loses stability via period doubling bifurcation, and a period 2 attractor appears so that the point $O$ becomes a saddle $O$ with $\dim W^u(O)=1$; and further
 %$\Rightarrow$
  the stable and unstable manifold of $O$ intersect (creating the ``homoclinic attractor''). Note also that the period-2 attractor must lose stability. In the case when the point $O$ has multipliers $\lambda_i$ with that $\lambda_1 < -1$, $0<\lambda_{2}<1, -1<\lambda_3<0$, $|\lambda_3|< \lambda_2$ and $|\lambda_1\lambda_2|>1$ the configuration of $W^u(O)$ is quite similar to those for the classical Lorenz attractor.\footnote{However, if $|\lambda_3|> \lambda_2$ and $|\lambda_1\lambda_3|>1$, then rather new subject, the so-called ``figure-eight'' pseudo-hyperbolic strange attractor can arise, see Fig.~\ref{Lor8attr}. Such discrete attractor was found in \cite{Kaz14} in a nonholonomic model of unbalanced ball moving on a plane.}

   This scenario provides us
with examples of truly high-dimensional (i.e. not two-dimensional) robust chaotic behavior which, as we sure, must occur in various applications. In the present paper we show that discrete Lorenz attractors can exist in  Poincar\'e map of a nonholonomic model of Celtic stone, see also \cite{Sasha13,GG13}.

\begin{figure}[htb]
\centerline{
  \includegraphics[width=12cm]{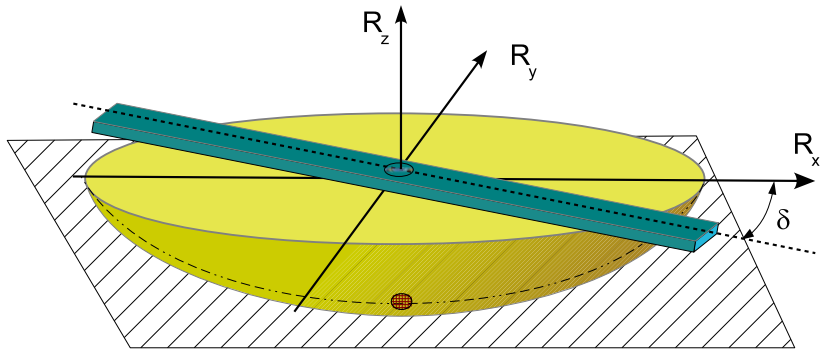}}
  \caption{
  {\footnotesize The sample of Celtic stone. The main body (of form of elliptic paraboloid) is symmetric with respect to the vertical axis $R_z$, the geometric horizontal axes are $R_x$ and $R_h$. The dynamical asymmetry is achieved due to a massive bar kept on the top of stone, which can be rotated by an angle $\delta$ with respect to the axis $R_x$. If $0<\delta<\pi/2$, the clockwise rotation of this stone is stable, whereas, the counterclockwise rotation is unstable.} }\label{realstone}
\end{figure}

Recall that, in the rigid body dynamics, the Celtic stone is a top (usually, of  symmetric form which remains a boat, see Fig.~\ref{realstone}) whose one of the principal inertial axes is vertical and other two axes are horizontal and  they are rotated by some angle with respect to the geometrical axes. A nonholonomic model of Celtic stone is a mathematical model provided that both the stone and the plane are absolutely rigid and rough, i.e. the stone moves along the plane without slipping and, moreover,  the friction force has zero moment. This means that the full energy is conserved which is a certain disadvantage of the model.
%dv oon som n absolutely rigid body which has
%
However, it is well known
that the nonholonomic model allows one to explain the main phenomenon of the Celtic stone dynamics -- the
nature of reverse, i.e., rotational asymmetry, which results in the fact that the stone
can rotate freely in one direction (e.g. clockwise) but ``does not want'' to rotate in the opposite direction (counterclockwise).
In the latter case it performs several rotations due to inertia, then stops rotating and starts oscillating,
after that it changes the direction of rotation and finally continues rotating freely (clockwise).

A mathematical explanation of this phenomenon seems now simple enough. The
fact is that, like most of the well-known nonholonomic mechanical models, the Celtic stone model
is described by a reversible system, i.e., a system that is invariant with respect to the coordinate and time change of the form  $X \to {\cal R}\; X, \;\; t \to -t$, where
${\cal R}$ is an involution, i.e. a specific diffeomorphism of the phase space such that ${\cal R}^2 = Id$.
However, in the case of Celtic stone,
this system is, in general, neither conservative nor integrable, although it possesses two independent integrals (see formula (\ref{eq:3-1})).
Because of this,
the system can possess, on the common level set of the integrals, asymptotic stable and completely
unstable solutions, stationary (equilibria), periodic (limit cycles) solutions etc., ${\cal R}$-symmetric with respect to
each other. Then, for example, a stable equilibrium corresponds to a stable vertical rotation of
the stone, and an unstable equilibrium symmetric with respect to it corresponds to an unstable
rotation in the opposite direction. Such an explanation of the reverse in the Celtic stone dynamics
was given in a series of
%recent
papers: by I.S.~Astapov \cite{A80}, A.V.~Karapetyan \cite{K81}, A.P. Markeev \cite{M83} and
others (further references can be found in
%the required
\cite{M2011}).

Nevertheless,  the motion of the Celtic stone
%on a plane
is still regarded in mechanics as one of the most complicated and poorly studied types of rigid body motion. Moreover, this is one of the few types of motion in which chaotic
dynamics was observed \cite{BM02, BM03,GGK12,GGK13}.\footnote{
Recently complex chaotic dynamics has been also found in a model of unbalanced rubber ball (a dynamically asymmetric ball with a displaced center of
gravity) moving on the plane, see \cite{K-RCD13}.}

The paper is organized as follows.
Section~\ref{sec:celt} contains necessary facts related to a nonholonomic model of Celtic stone under consideration.
Section~\ref{sec:2stone} is the main part of the paper. In this section we consider the nonholonomic model of Celtic stone with parameters (\ref{eq:fisik2}) and show (mostly numerically) that a discrete Lorenz attractor exists in this model. We study this problem using a research strategy (Section~\ref{sec:32})
%containing four main steps)
based on
%the definition (see Definition~1) of a pseudo-hyperbolic diffeomorphism and
fundamental facts from the theory of Lorenz-like attractors (see Section~\ref{sec:31}). Main results are formulated and discussed in Section~\ref{sec:33}.

\section{The nonholonomic model of Celtic stone.} \label{sec:celt}

We study the dynamics of a rigid body moving on a plane without slipping.
This means that we consider a nonholonomic model of motion in which the
contact point of the body has zero velocity, i.e. we have
\begin{equation}
 {\vv} %\vv
+ \vo \times \vr = 0
\label{eq:1cs}
\end{equation}
where $\vr$ is the radius vector from the center
of mass $C$
to the contact point,
$\vv$ is the velocity of $C$
%the center of mass
and $\vo$ is the angular velocity of the body. As usual, the coordinates
of all vectors are defined in some coordinates rigidly attached to the
body. Then the equations of motion can be written in the form
\cite{BM01-2}
%
%\iffalse
\begin{equation}
\begin{array}{l}
\dot{\vM} =
\vM\times\vo +
m\dot{\vr}\times (\vo\times\vr) + m{\mbox g}\vr\times\vg, \\
\dot{\vg} = \vg\times\vo,
\end{array}
\label{eq:2}
\end{equation}
%\fi
where
\begin{equation}
\label{eq:2-1}
\vM = \vI\;\vo + m{\vr}\times (\vo\times\vr)
\end{equation}
is the angular momentum of the body with respect to the contact point,
$\vg$ is the unit vertical vector and $m {\mbox g}$ is the gravity force.
The equation (\ref{eq:2}) admits two integrals
%(the {\em full mechanical energy} and the {\em geometric} integrals correspondingly):
%
\begin{equation}
{\cal H} = \frac{1}{2}({\vM}, \vo) - m{\mbox g}({\vr},\vg)\;\; \mbox{and} \;\; (\vg, \vg) = 1,
\label{eq:3-1}
\end{equation}
{\em the energy integral}
%(the first integral)
and {\em the geometric integral}, respectively.

We consider the Celtic stone whose surface $F(\vr)$ has the shape of the {\em
elliptic paraboloid}
$$
\displaystyle F(\vr) = \frac{1}{2}\left(\frac{r_1^2}{a_1} + \frac{r_2^2}{a_2}\right) - (r_3 + h) = 0,
$$
where $a_1$ and $a_2$ are the principal radii of curvature at the
paraboloid vertex $(0,0,-h)$ respectively, and the center of mass is the point $r_1 =
r_2 = r_3 = 0$.
Therefore, the vector $\vr$ and $\vg$ are related by:
\begin{equation}
\begin{array}{l}
\displaystyle r_1 = - a_1\frac{\gamma_1}{\gamma_3}, \;\; r_2 = -a_2\frac{\gamma_2}{\gamma_3},\;\;
r_3 = -h + \frac{a_1\gamma_1^2 + a_2\gamma_2^2}{2\gamma_3^2}.
\end{array}
\label{eq:5}
\end{equation}

It is also assumed that
%the center of mass is the point $r_1 = r_2 = r_3 = 0$ and
one of the principal axes of inertia is vertical.
One of the main features of Celtic stone is that two other principal
axes of inertia are rotated about the geometrical axes by some angle
$\delta$, where $0 < \delta < \pi / 2$.
Accordingly, the inertia tensor takes the following form \cite{BM03}:
\begin{equation}
\vI \;=\; \left(
\begin{array}{ccc}
I_1 \cos^2\delta + I_2\sin^2\delta & (I_1 - I_2) \cos\delta \sin\delta & 0 \\
(I_1 - I_2)\cos\delta \sin\delta & I_1\sin^2  \delta + I_2 \cos^2 \delta &  0 \\
0 & 0 & I_3
\end{array}
\right)\;\;,
\label{eq:6}
\end{equation}
where $I_1, I_2$ and $I_3$ are the principal moments of inertia of the
stone.
%\iffalse
We express vectors $\vr$, $\dot{\vr}$ and $\vo$ by $\vM$ and $\vg$ using relations  (\ref{eq:1cs}), (\ref{eq:2-1}), (\ref{eq:5}) and (\ref{eq:6}).
Then the system (\ref{eq:2}) can be represented in the standard form
\begin{equation}
(\dot{\vM}, \dot{\vg}) = G(\vM, \vg, \mu),
\label{eq:1-1}
\end{equation}
%
%\fi
of six-dimensional system with respect to phase variables $\vM$ and $\vg$.
%\in \mathbb{R}^3$ and $\vg \in \mathbb{R}^3$:
%
%
This system depends also on parameters $\mu$ characterizing the geometrical
and physical properties of the stone.
Note that on the common level set of the integrals (\ref{eq:3-1}) the system (\ref{eq:2})
defines the flow on a four-dimension manifold:
$
{\cal M}^4 = \{(\vM, \vg): (\vg, \vg) = 1,
{\cal H}(\vM, \vg) = \mbox{const}\},
$
which is homeomorphic to ${\mathbb{S}^3 \times \mathbb{S}^3}$.

\subsection{The Andoyer-Deprit variables.} \label{sec:AD}

In numerical investigations of dynamics of Celtic stone we use the
so-called {\em Andoyer-Deprit} variables $(L, H, G, g, l)$ defined by the
formulas \cite{BM01-2}

\begin{equation}
\begin{array}{l}
\displaystyle M_1 = \sqrt{G^2 - L^2}\sin l, \;\;M_2 = \sqrt{G^2 - L^2}\cos l, \;\;M_3 = L, \\ %\\
%
%\displaystyle \gamma_1 = \left(\sin\zeta\cos\tau + \sin\tau\cos\zeta\cos g\right)\sin l + \cos\zeta\sin g\cos l, \\
\displaystyle \gamma_1 = \left(\frac{H}{G} \sqrt{1-\frac{L^2}{G^2}} + \frac{L}{G} \sqrt{1-\frac{H^2}{G^2}} \cos g\right) \sin l + \sqrt{1-\frac{H^2}{G^2}} \sin g \cos l \\
\displaystyle \gamma_2 = \left(\frac{H}{G} \sqrt{1-\frac{L^2}{G^2}} + \frac{L}{G} \sqrt{1-\frac{H^2}{G^2}} \cos g\right) \cos l - \sqrt{1-\frac{H^2}{G^2}} \sin g \sin l \\
\displaystyle \gamma_3 = \frac{HL}{G^2} - \sqrt{1-\frac{L^2}{G^2}} \sqrt{1-\frac{H^2}{G^2}} \cos g
%\left(\frac{H}{G}\sqrt{1-\left(\frac{L}{G}\right)^2} +  \frac{L}{G}\sqrt{1-\left(\frac{H}{G}\right)^2}\cos g \right) \sin l +
%
%
\end{array}
\label{eq:6-AD}
\end{equation}
%where $\displaystyle  \sin\tau = {L}/{G},\;\; \sin\zeta = {H}/{G}$,
By definition (see, e.g., \cite{BM01-2}),
\begin{equation}
H = (\vM, \vg) = M_1\gamma_1 + M_2\gamma_2 + M_3\gamma_3.
\label{eq:6-HM}
\end{equation}

On the common level set of two integrals (\ref{eq:3-1}),
the system (\ref{eq:1-1}) represents the four-dimensional flow ${\cal
G}_{E}$. Note that the new coordinates $L, H, G, g$ and $l$ are chosen in
such a way that the condition $(\vg, \vg) = 1$ holds automatically. Thus,
the formulae (\ref{eq:6-AD}) specify a one-to-one correspondence between
the coordinates $\{ (\vM, \vg) : \vg^2 = 1\}$ and $(L, H, G, g, l)$
%$\{\boldsymbol{M}, \boldsymbol{\gamma}: \boldsymbol{\gamma^2} = 1\}$
everywhere except for the planes
%\begin{equation}
${L}/{G} = \pm 1$ and ${H}/{G} = \pm 1$ (for which the coordinate $l$ and, respectively, $g$ are not defined).

Further, we will investigate the systems on the four-dimensional energy
levels\\ ${\cal H}(L, G, H, l, g) = E$. In this case the planes $g = g_0 =
\mbox{const}$ (for appropriate $g_0$) can be considered as cross-sections
for orbits of the corresponding four-dimensional flow ${\cal G}_{E}$.
Thus, we can also study the dynamics of the three-dimensional Poincar\'e map
\cite{BM02, BM03}:
\begin{equation}
\bar x = {\cal F}_{g_0}(x), \;\; x = \left(l,\frac{L}{G},\frac{H}{G}\right),
\label{eq:mapPoin}
\end{equation}
which is defined in the domain $0 \leq l < 2\pi, -1 < \frac{L}{G} < 1, -1
< \frac{H}{G} < 1$.

\subsection{Symmetries in the Celtic stone model.}

The system (\ref{eq:1-1}) possesses a number of interesting and useful
symmetries described by the following lemma.

\begin{lm} {\rm \cite{BM03}}
In the case under consideration, the system (\ref{eq:1-1}) is symmetric with respect
to the coordinate changes:
\begin{equation}
%
%\qquad
\hspace{-1cm}  (a) \qquad {\cal S}_1:\; \omega \to -(-\omega_1,-\omega_2,\omega_3)\;\;,\;\;\gamma\to (- \gamma_1,- \gamma_2,\gamma_3) \\
\label{eq:inv-sym}
\end{equation}
and is reversible with respect to the following involutions:
\begin{equation}
\begin{array}{l}
(b)\qquad  {\cal I}_1:\; \omega \to - \omega, \;\;,\;\;\gamma\to  \gamma,\;\;,\;\; t\to - t\\
%\label{eq:inv-I1}
%
%
%
(c)\qquad {\cal I}_2:\; \omega \to (\omega_1,\omega_2,-\omega_3)\;\;,\;\;\gamma\to (- \gamma_1,- \gamma_2,\gamma_3),\;\; t\to -t
\end{array}
\label{eq:inv-inv}
\end{equation}

\label{lm:sym1}
\end{lm}

Note that these symmetries and involutions are also preserved for the
Andoyer-Deprit coordinates. However, they cannot always be linear in this
case. But for the Poincar\'e map (\ref{eq:mapPoin}) with the
crosssection $g = 0$, which we denote as ${\cal F}_0$, the symmetries
(\ref{eq:inv-inv}) remain linear.
\begin{lm}  {\rm \cite{BM03}}
The map ${\cal F}_0$ is invariant under the following transformations:
\begin{equation}
\begin{array}{l}
\smallskip %vspace{5mm}
\displaystyle \hspace{-0.0cm} (a) \qquad\;\;  \tilde{\cal S}_1 :\;\; l \to l + \pi,\;\;\frac{L}{G} \to  \frac{L}{G},\;\;\frac{H}{G}\to  \frac{H}{G}, \\
\smallskip %vspace{5mm}
\displaystyle (b)\;\; \tilde{\cal I}_1:\; l \to l + \pi,\;\;\frac{L}{G}\to - \frac{L}{G},\;\;\frac{H}{G}\to - \frac{H}{G},\;\; {\cal F}_0\to {\cal F}_0^{-1} \\
%
%\vspace{5mm}
%
\displaystyle (c)\;\; \tilde{\cal I}_2 = \tilde{\cal I}_1\tilde{\cal S}_1  :\; l \to l,\;\;\frac{L}{G}\to - \frac{L}{G},\;\;\frac{H}{G}\to - \frac{H}{G},\;\; {\cal F}_0\to {\cal F}_0^{-1}
\end{array}
\label{eq:inv}
\end{equation}
\label{lm:sym2}
\end{lm}

\begin{crl}
Let $L^*$ be an orbit of ${\cal F}_0$. Then $\tilde S_1(L^*)$, $\tilde{\cal
I}_1(L^*)$ and $\tilde{\cal I}_2(L^*)$ are also orbits of ${\cal F}_0$.
Moreover, the orbits $L^*$ and $\tilde S_1(L^*)$, as well as $\tilde{\cal
I}_1(L^*)$ and $\tilde{\cal I}_2(L^*)$, are symmetric with respect to each
other. The orbits $L^*$ and $\tilde S_1(L^*)$ are both in involution with
the orbits $\tilde{\cal I}_1(L^*)$ and $\tilde{\cal I}_2(L^*)$.
\label{cor:SyInv}
\end{crl}

\section{On discrete Lorenz attractors in Celtic stone dynamics.} \label{sec:2stone}

In this section we consider the nonholonomic model of Celtic stone whose physical parameters are
as follows:
 \begin{equation}
I_1 =2, I_2 = 6, I_3 =7, m=1, {\rm g} = 100, a_1 = 9, a_2=4,  h=1.
\label{eq:fisik2}
\end{equation}
We also take  $\delta = 0.485$.

Note that the Celtic stone model with the parameters (\ref{eq:fisik2}) was considered in \cite{KJSS12} in which a strange attractor was found, with $E=770, \delta=0.405$. This strange attractor is quite similar to the
attractor in Fig.~\ref{fig3-3}(f).
Since analogous attractors are known to exist in three-dimensional H\'enon
maps \cite{GGS12,GGOT13} near the  boundaries of destruction of discrete Lorenz attractors,
the question naturally arises whether a discrete Lorenz attractor exists for close values of the parameters $E$ and $\delta$.
The answer is positive and we give below a short review of results
obtained.

Before description of main results (Sections~\ref{sec:32} and~\ref{sec:33})we make some remarks and recall necessary definitions.

\subsection{On pseudo-hyperbolic Lorenz-like attractors.}  \label{sec:31}

First of all
we give the definition of pseudo-hyperbolicity for diffeomorphisms, which is, in fact,
a reformulation  of the corresponding definition for flows from \cite{TS08}, see also \cite{Sat}.

Let $f$ be a $C^r$-diffeomorphism, $r\geq 1$ and let $Df$ be its linearization.
%: $Df(x_0) = \left(\frac{\partial f}{\partial x}_{|x=x_0}\right)$.
An open bounded
domain ${\cal D} \subset \mathbb{R}^n$ is \textsf{absorbing} for
$f$ if $f(\overline{{\cal D}}) \subset {\cal D}$.

\begin{df}
The diffeomorphism $f$ is called \textsf{pseudo-hyperbolic} on $D$ if the following conditions
hold.
\begin{itemize}
\item[{\rm 1)}] For each point of $D$ there exist two  transversal
subspaces $N_1$ and $N_2$ continuously depending on the point
($\dim N_1 = k\geq 1, \dim N_2 = n-k$) which are invariant with respect to $Df$:
    $$
    Df(N_1(x)) = N_1(f(x)),\;\;Df(N_2(x)) = N_2(f(x)).
    $$
and for each orbit $L:\{x_i \;|\; x_{i+1}= f(x_i),\; i=0,1,...; x_0\in D \}$ the maximal
Lyapunov exponent corresponding to the subspace $N_1$ is strictly smaller than the minimal Lyapunov exponent
corresponding to the subspace $N_2$, i.e., the following inequality holds:

\begin{equation}
\begin{array}{l}
\displaystyle \lim\limits_{n\to\infty}\sup \frac{1}{n}\;\ln\; (\!\!\sup\limits_{\begin{array}{c}u\in N_1(x_0)\\ \|u\|=1\end{array}}\!\! \|D^nf(x_0) u\|) < \\
\qquad\qquad\qquad
\displaystyle <\;\lim\limits_{n\to\infty}\inf \frac{1}{n}\;\ln \;(\!\!\inf\limits_{\begin{array}{c} v\in N_2(x_0)\\ \|v\|=1\end{array}}\!\! \|D^nf(x_0) v\|).
\end{array}
\label{eq:psh1}
\end{equation}

\item[{\rm 2)}]  The restriction of $f$ to $N_1$ is contracting, i.e., there exist
constants $\lambda>0$ and $C_1>0$ such that
\begin{equation}
\|D^n f(N_1)\| \leq  C_1 e^{-\lambda n}.
\label{eq:psh2}
\end{equation}

\item[{\rm 3})]  The restriction of $f$ to $N_2$ expands volumes exponentially, i.e., there exist
such constants $\sigma>0$ and  $C_2>0$ such that
\begin{equation}
|\det D^n f(N_2) | \geq  C_2 e^{\sigma n}.
\label{eq:psh3}
\end{equation}
\end{itemize}
\label{def:psevhyp}
\end{df}

The following property immediately follows from this definition:

\begin{itemize}
\item[$1^*$]
All the orbits in ${\cal D}$ are unstable: each orbit has positive maximal Lyapunov
exponent
$$
\Lambda_{max}(x) = \lim\limits_{n\to\infty}\sup \frac{1}{n}\;\ln\;  \|D^nf(x)\| > 0
$$
\end{itemize}

Note that the pseudo-hyperbolicity conditions require the expansion of only $(n-k)$-dimensional volumes by the restriction of the diffeomorphism to $N_2$,
which makes these conditions different from those for uniform hyperbolicity, where
the following condition must hold: $\|D^{-n}f(N_2)\| < C e^{-\sigma n}$, i.e., the uniform
expansion should be along all directions in $N_2$. Nevertheless,  it is possible to
establish the following fact in a standard way \cite{An67,TS98}.

\begin{itemize}
\item[$2^*$]
The pseudo-hyperbolicity conditions are not violated under small $C^r$-perturbations
of the system. Moreover, the spaces $N_1$ and $N_2$ change continuously.
\end{itemize}

These two conditions imply that if the diffeomorphism $f$ has an attractor in ${\cal D}$, then
this attractor is strange and does not contain stable periodic orbits, which also do not appear
under small perturbations. In other words, pseudo-hyperbolic attractors are {\em genuine strange attractor}.
The discrete Lorenz attractors
%Lorenz-like attractors for diffeomorphisms, we call them also ,
form a certain subclass of the class of
pseudo-hyperbolic attractors.

Note that the dynamical properties of the geometric Lorenz model \cite{ABS83}
under small time-periodic perturbations were investigated in \cite{TS08}.
It was also shown that the properties of pseudo-hyperbolicity and
chain transitiveness\footnote{i.e. when any two points in an invariant set can be joined by an $\varepsilon$-orbit belonging to the set, for any sufficiently small $\varepsilon$} of a
non perturbed Lorenz attractor hold for a periodically perturbed attractor as well.
Thus, the Poincar\'e map
%(the map for a period of perturbation)
also possesses here
a pseudo-hyperbolic attractor $A$, which appears to be a basic example of a {\em discrete Lorenz attractor}.
Note that the saddle equilibrium of the Lorenz attractor corresponds, after perturbations, to
%  of the unperturbed flow corresponds now to
the saddle type fixed point of the corresponding Poincar\'e map.

The same conclusions can also be drawn without assuming that the map under consideration is
a Poincar\'e map of a system periodic in time and close to an autonomous one. For this general case
the corresponding definition of  a discrete Lorenz attractor was given in
\cite{GGOT13}. Formally, this definition requires to satisfy main conditions which hold for the Poincar\'e map
of periodically perturbed Lorenz attractor, including the geometrical analogy and the spectrum of Lyapunov exponents due to the Definition~\ref{def:psevhyp}.
Thus, for numerical study of a discrete Lorenz attractor, we need to check the required geometrical analogy and that fact that the {\em numerical} Lyapunov exponents
$\Lambda_1,\Lambda_2,\Lambda_3$ satisfy conditions  $\Lambda_3 < \min\{\Lambda_1,\Lambda_2\}$, $\Lambda_3 < 0$ and $\Lambda_1 + \Lambda_2 >0$
(some analogs of conditions (\ref{eq:psh1}), (\ref{eq:psh2}) and (\ref{eq:psh3}), respectively).

In Section~\ref{sec:32} we describe the corresponding research strategy.The obtained results are collected in Section~\ref{sec:33}.

\subsection{A strategy of qualitative and numerical
study  of discrete Lorenz attractors on example of the family ${\cal F}_{0E}$.} \label{sec:32}

When studying dynamics and bifurcations in the family ${\cal F}_{0E}$ of the Poincar\'e map ${\cal F}_{0E}$
(\ref{eq:mapPoin}) for appropriate values of $E$ (with fixed $\delta = 0.485$), we will act by employing the following strategy which, in fact, is justified by
%the given above
Definition~\ref{def:psevhyp} and its corollaries $1^{*}$ and $2^{*}$.

\begin{itemize}
\item[1)]
Verify the geometrical similarity of our attractor $A_{E^*}$ found in the Celtic stone model
to the Lorenz attractor. Here the strange attractor which was found for $E=E^*=752$
is examined.
\end{itemize}

\begin{figure}[ht]
\begin{tabular}{cc}
\psfig{file=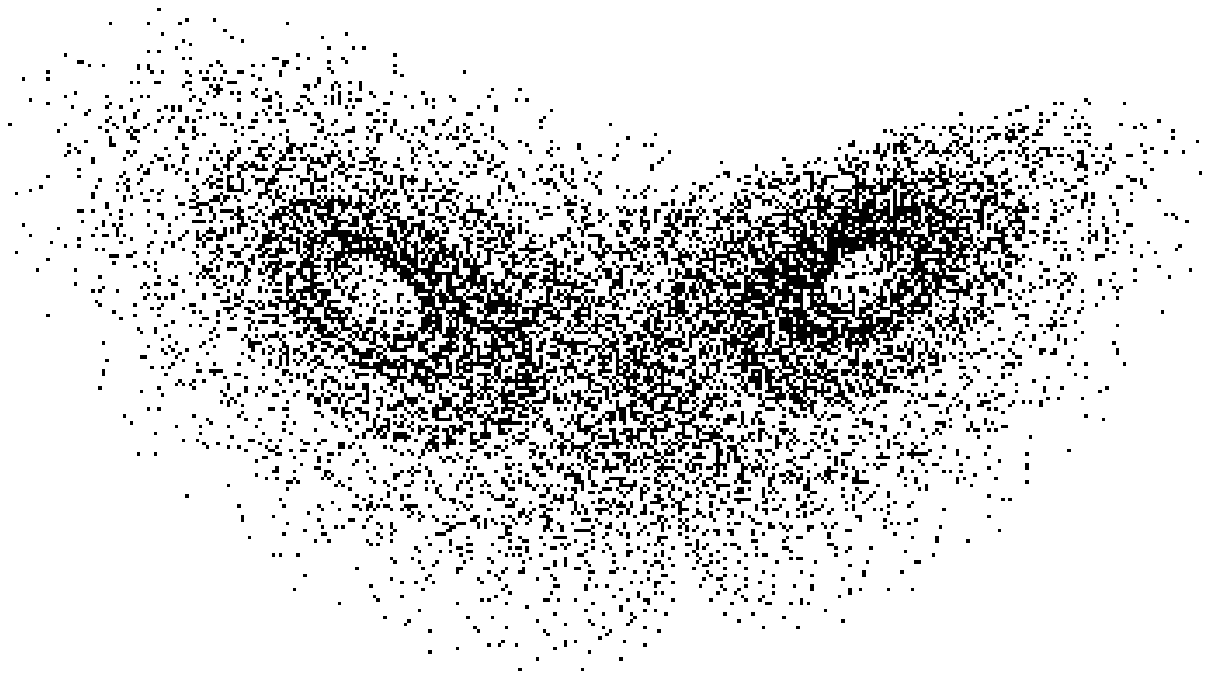,width=75mm} &
%
%\hspace{15mm}
\psfig{file=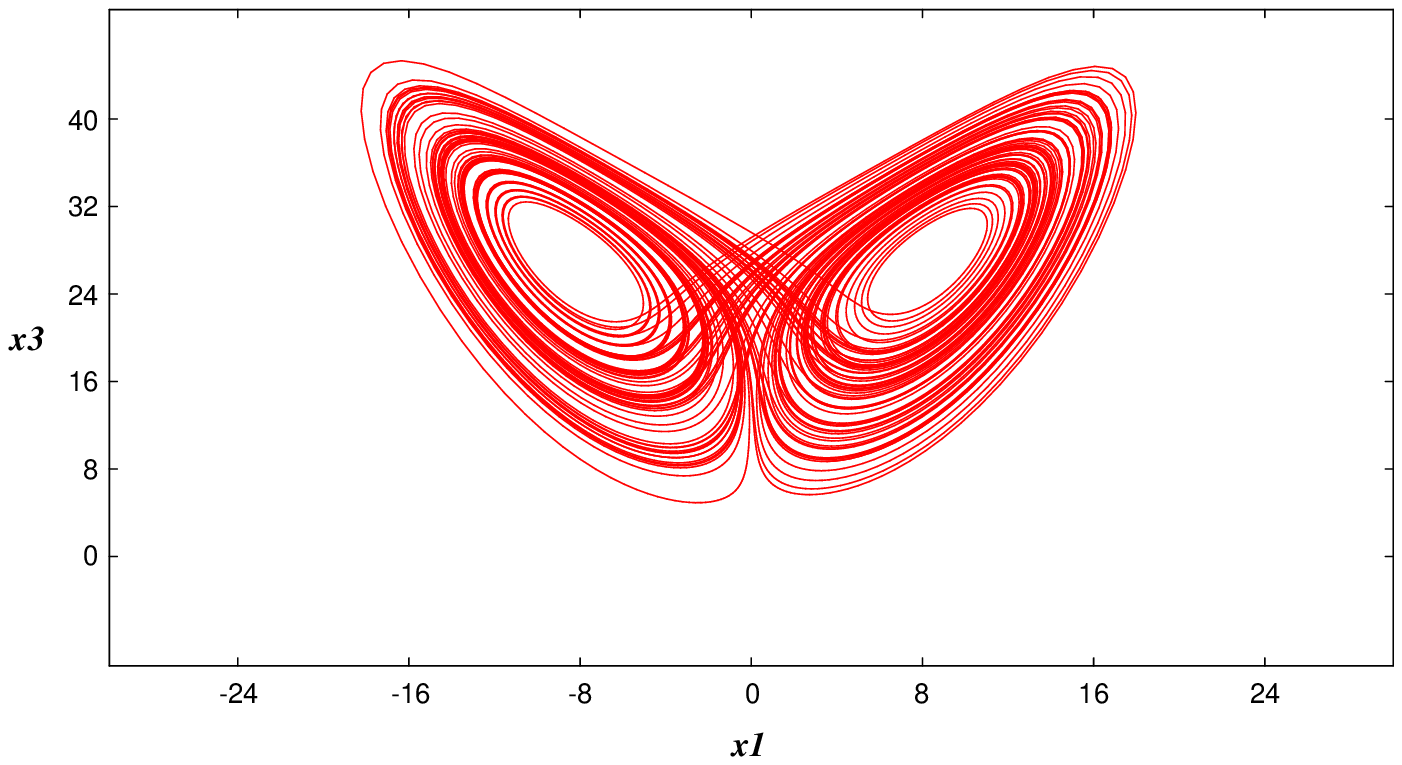, width=70mm} \\
(a)&(b)
\end{tabular}
\caption{{\footnotesize  (a)  a Lorenz-like attractor for $E=E^*=752$ in the Celtic stone model
(about 10000 iterations of some initial point are shown);  b)  the projection of the Lorenz
attractor from the Lorenz model onto the $(x,z)$ plane.}
}
\label{fig3-4}
\end{figure}

In particular, this similarity manifests itself in the fact that our three-dimensional map ${\cal F}_{0E^*}$ possesses the following
features: (i) it has a fixed saddle point
$O^*$ belonging to the attractor $A_{E^*}$ with the multipliers of $\lambda_1, \lambda_2, \gamma$ such that
$|\lambda_2|< |\lambda_1| < 1 < |\gamma|$,  $\lambda_1 > 0, \lambda_2 < 0, \gamma < -1$ and
$|\lambda_1 \gamma| > 1$; (ii) the manifolds $W^u(O)$ and $W^s(O)$ have nonempty intersection;
(iii) the phase portraits look ``similar'', see Fig.~\ref{fig3-4}, and, moreover, the main bifurcations leading to appearance of the attractor are quite analogous to those in the Lorenz system, see item 4) below.

We mention that negative values of the multipliers  $\lambda_2$ and $\gamma$ provide
the Lorenz symmetry ($x \to x, y \to -y, z \to -z)$ of the homoclinic structure.
Moreover, for the values of the parameter $E$ close to $E^*$, the manifold $W^u$ will intersect
$W^s$ strictly from one side of the strong stable invariant manifold $W^{ss}(O)$,
which is tangent to the eigendirection corresponding to the multiplier $\lambda_2$ of $O^*$,
which provides the
%обеспечивает
homoclinic configuration of the figure-eight-butterfly similar to the Lorenz attractor.

\begin{itemize}
\item[2)]
Verify numerically the strangeness and pseudo-hyperbolicity of the attractor $A_{E^*}$.
%(we refer to Section~\ref{sec:appendix} for the definition of pseudo-hyperbolicity and related terms).
\end{itemize}

At this stage we investigate the spectrum  $\Lambda_1,\Lambda_2,\Lambda_3$ of the Lyapunov
exponents of the map ${\cal F}_{0E^*}$ on the attractor $A^*$ and show that this spectrum, where
$\Lambda_1 >\Lambda_2 >\Lambda_3$, satisfies the following conditions: (1) $\Lambda_1 > 0$;
(2) $\Lambda_1 +\Lambda_2 + \Lambda_3<0$; (3) $\Lambda_1 + \Lambda_2 >0 $.
The conditions (1) and (2) imply that the attractor $A^*$ is strange and the condition (3) holds when the attractor
%implies that it
is pseudo-hyperbolic (the map expands two-dimensional areas transversal to the strong
contraction direction related to the exponent $\Lambda_3 < 0$).

\begin{itemize}
\item[3)]
Plot numerically the dependence of the maximal exponent $\Lambda_1$ for some range of
the parameter $E$ containing this value, $E = E^*$, for which the attractor $A_{E^*}$ exists, see
Fig.~\ref{WU-Lorkelt}.
\end{itemize}

At this stage we verify (only numerically) that our attractor is not a quasi-attractor, i.e., it
does not contain stable periodic orbits of large periods, which do not appear under perturbations either.
As is seen from Fig.~\ref{WU-Lorkelt}, the graph resides in the domain $\Lambda_1>0$
and looks like a continuous function, whereas, if  $A_{E^*}$ were a quasi-attractor,
the ``holes'' would be observed on the plot containing the ranges of $\Lambda_1<0$
corresponding to ``stability windows''.

\begin{itemize}
\item[4)]
Investigate (mostly numerically) the main bifurcations starting at the stable fixed point and
leading to the appearance of discrete Lorenz attractors (including  $A_{E^*}$). We also trace the main stages of destruction of strange attractor.
\end{itemize}

In principle, this item may seem unnecessary but we suppose it to be the most interesting
because here one can follow a certain ``genetic'' connection between the phenomena observed
in flows with the Lorenz attractors (Lorenz model, Shimizu-Morioka model etc.) and those observed in the
model of Celtic stone.
Moreover, as the calculations show, our Poincar\'e map ${\cal F}_{0E}$ behaves like a ``small perturbation of
the time shift of the flow in the geometric Lorenz model'' for corresponding values of $E$.
Formally, this circumstance can be caused by the interesting fact that the middle Lyapunov exponent $\Lambda_2$
is very close to zero (for the flow case it is simply equal to zero):
during calculations it demonstrates small oscillations in the range between $0.00007$ and
$0.00015$, see \cite{GOST05} for a discussion of this topic. But what is really interesting is that
the bifurcations leading to the appearance of strange attractor are here almost identical to those which accompany
the birth of strange attractor in the Lorenz model \cite{Sh80},
see.
%\S~\ref{sec:chis} и
Fig.~\ref{fig3-2a}(b)--(f).

\subsection{Results of the numerical study.}  \label{sec:33}

Below we show the results of numerical investigations performed according to items
1)--4) of the strategy.

 \begin{figure}[ht]
\centerline{\epsfig{file=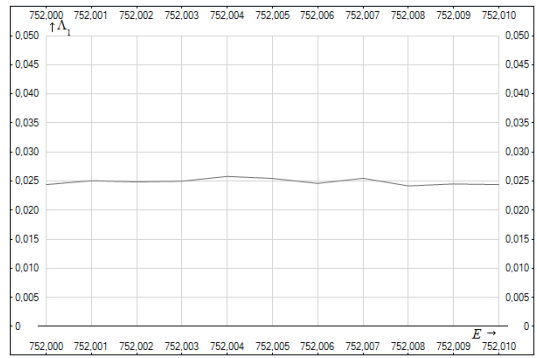, width=7cm
%, height=10cm
}} \caption{
{\footnotesize The plot of $\Lambda_1(E)$ in the range of $[752;752.01]$.
}}
\label{WU-Lorkelt}
\end{figure}

1) Fig.~\ref{fig3-4} shows  (a) iterations of a single point of the attractor $A_{E^*}$
of the map $T_E$ for $E =E^* =752$ (for an appropriate angle of projection) and (b) the projection of
orbits of the classical Lorenz attractor from the Lorenz model for $r = 28$, $\sigma = 10$, and $b = 8/3$
onto the $(x,z)$ plane displayed for comparison.
%Очевидно, что здесь прослеживается определенное сходство.
%\newpage

The fixed saddle point
$O^*$ with coordinates of $l=3.650;\; {L}/{G} = 0.669;\;{H}/{G} = -0.384$ on the
attractor $A_{E^*}$ has the multipliers
$\lambda_1 = 0.996;\; \lambda_2 = -0.664;\; \gamma = -1.312$.
If one draws its unstable manifolds (``separatrices''), then, as expected, they will have
``loops'' (due to the existence of the homoclinic intersection), see Fig.~\ref{fig3-2b}(a),
in contrast with the unstable separatrices of the Lorenz attractor in flows which appear to be
sufficiently monotonous spirals.

2) For the attractor $A_{E^*}$ at $E= E^*=752$ the spectrum of the Lyapunov exponents
was obtained as follows:
$\Lambda_1 = 0.0248;\;\; \Lambda_3 = - 0.2445,\;\; 0.00007<\Lambda_2<0.00015$.

Evidently, the conditions $\Lambda_1>0$,  $\Lambda_1 +\Lambda_2+\Lambda_3<0$ and  $\Lambda_1 + \Lambda_2 >0$ hold here.

3) On the graph of Fig.~\ref{WU-Lorkelt} the dependence of the maximal Lyapunov exponent $\Lambda_1 = \Lambda_1(E)$ on $E$
is shown for the range $[752;752.01]$ of the parameter $E$.

\begin{figure}[htb]
\begin{tabular}{ccc}
\psfig{file=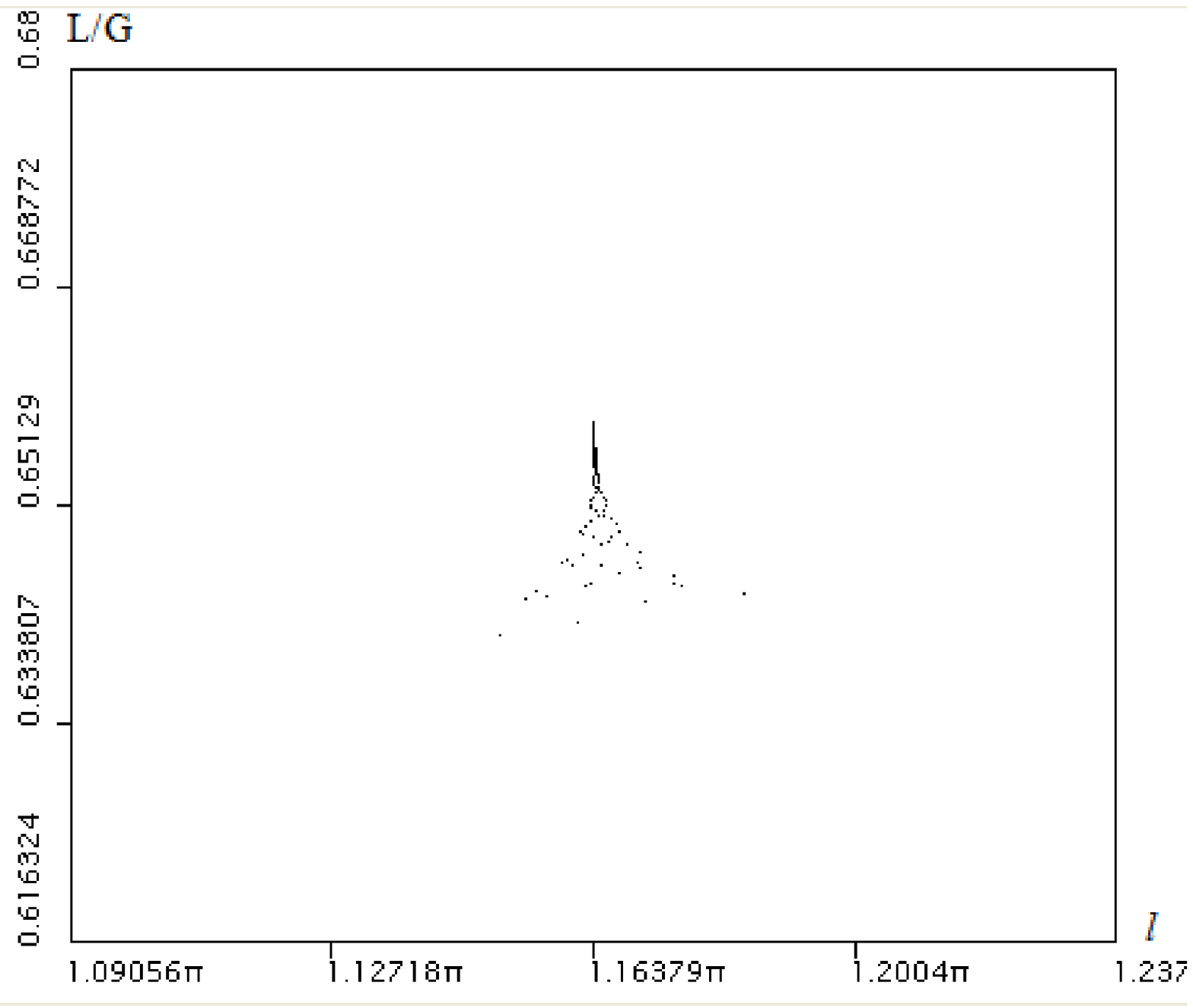,width=50mm} &
%
%\hspace{-5mm}
\psfig{file=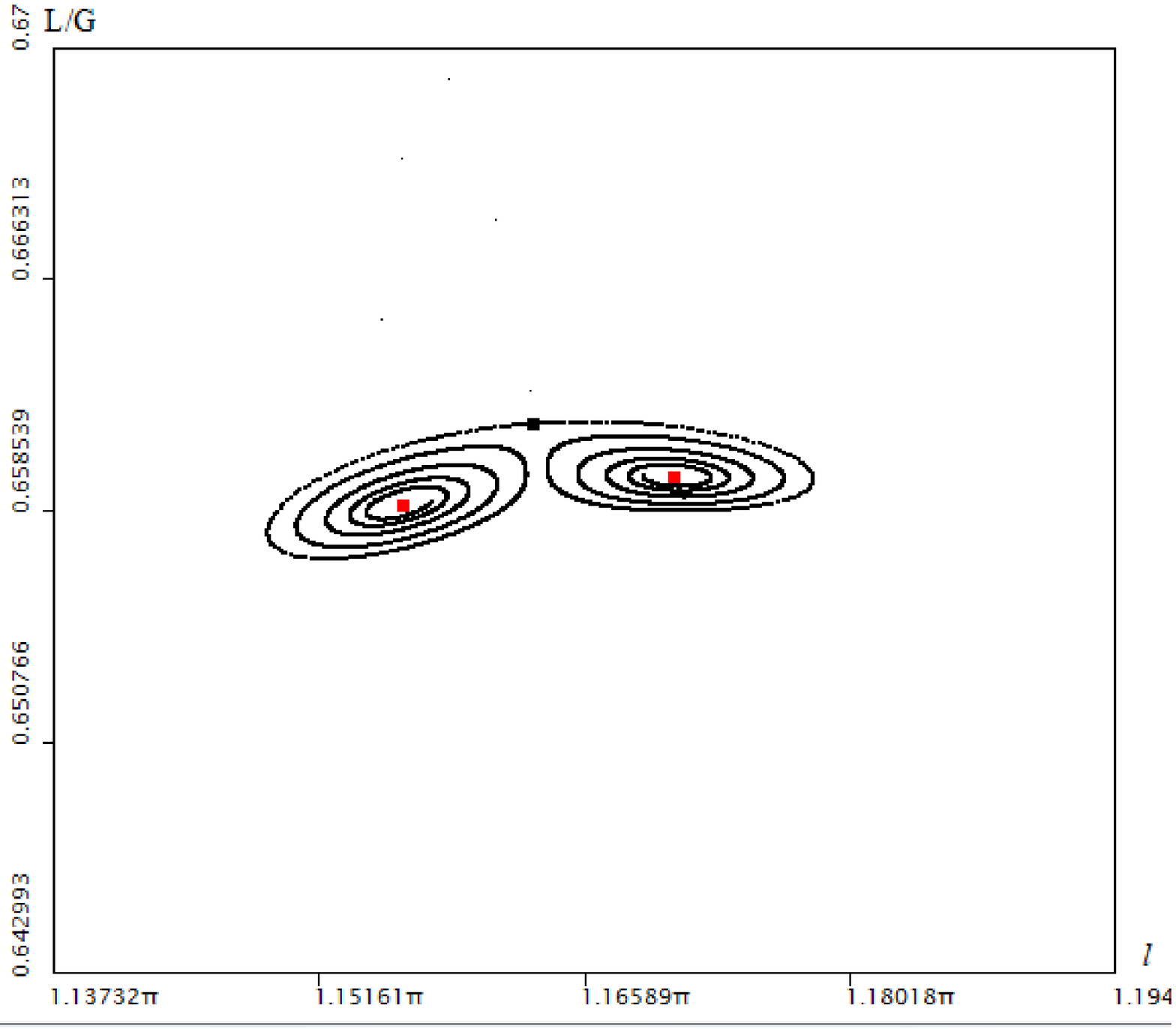, width=50mm}&
\psfig{file=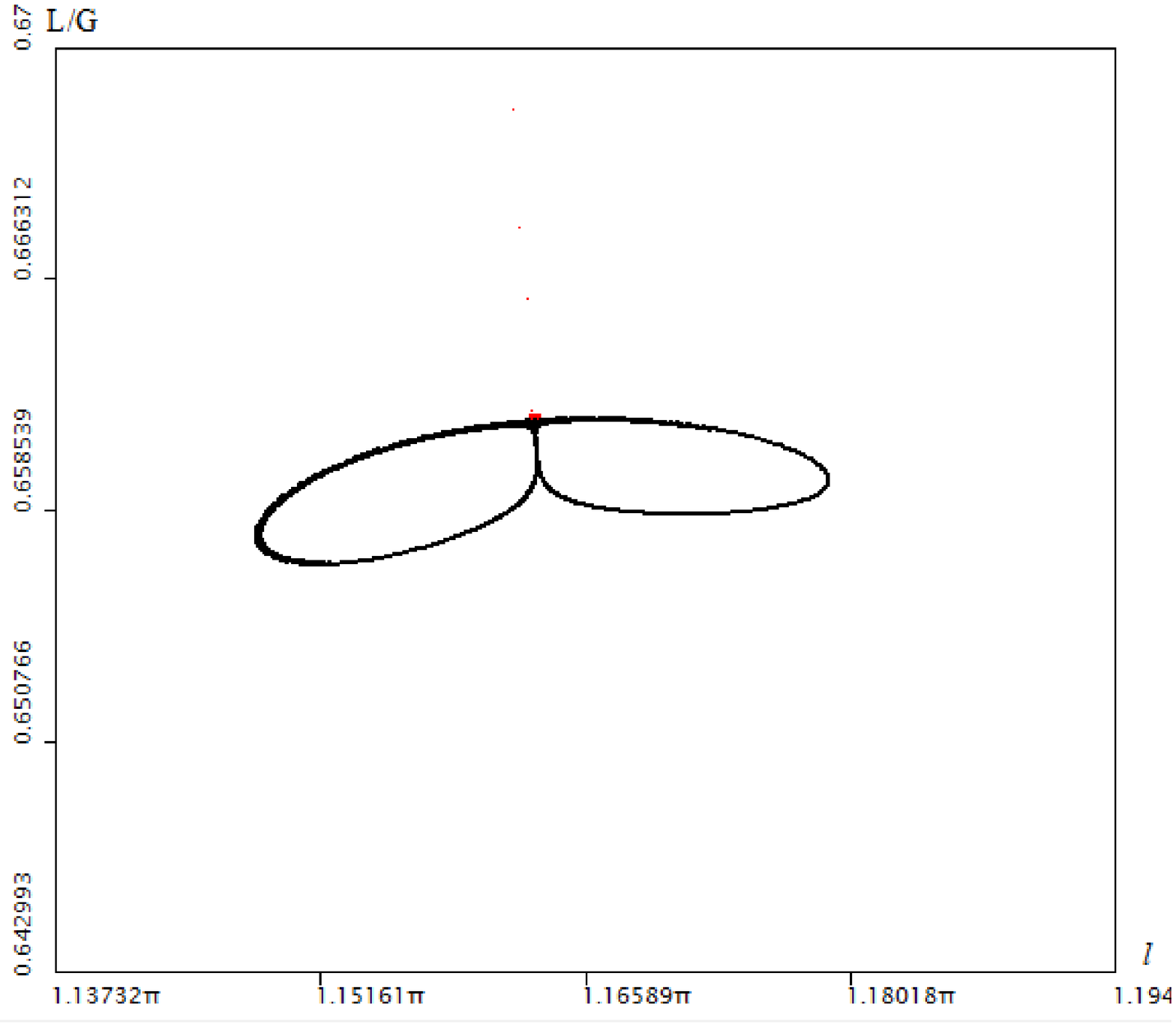,width=50mm} \\
a) E=747 &  b) E=748.4 & c) E=748.4395 \\
\psfig{file=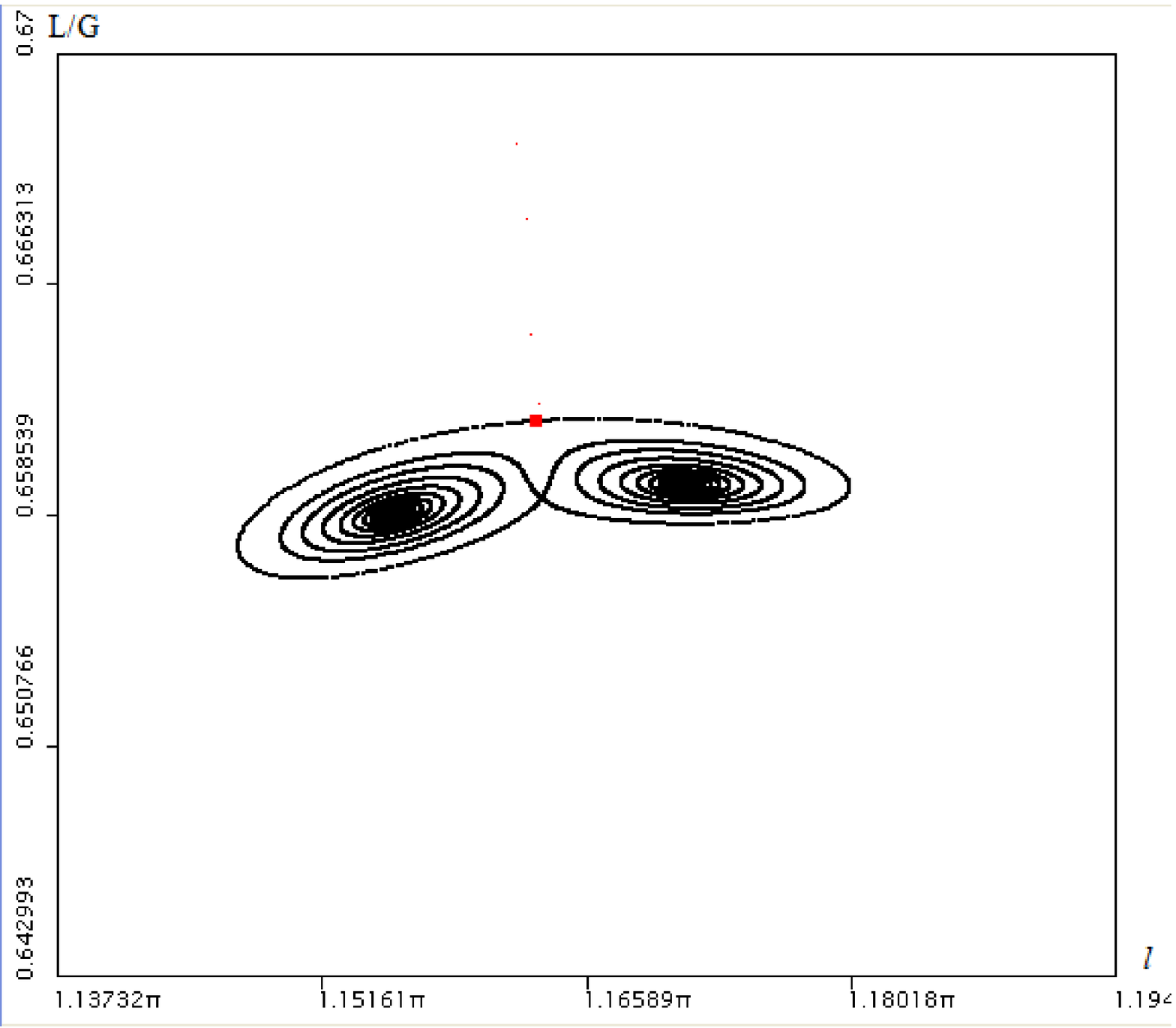, width=50mm}&
\psfig{file=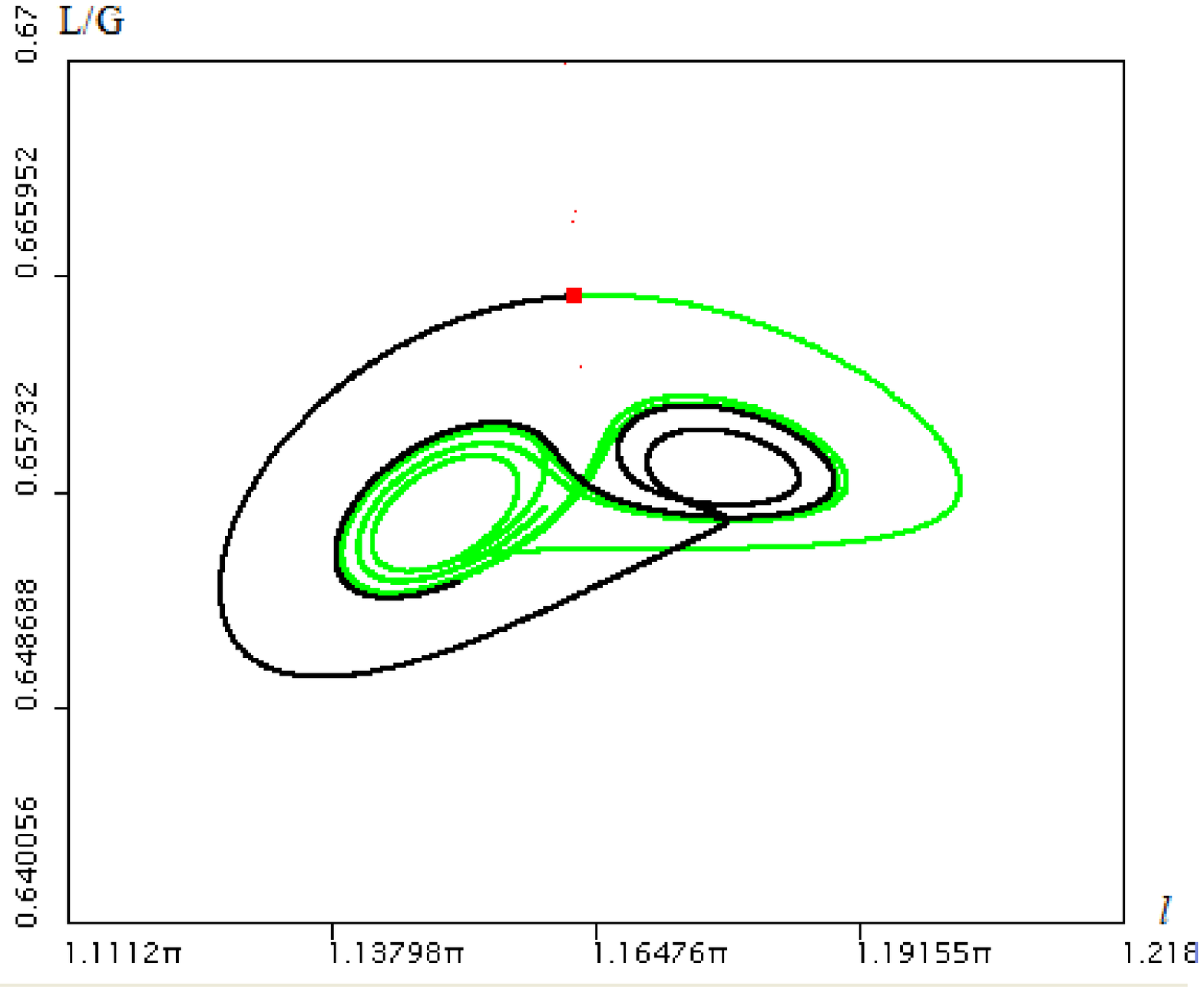,width=50mm}&
\psfig{file=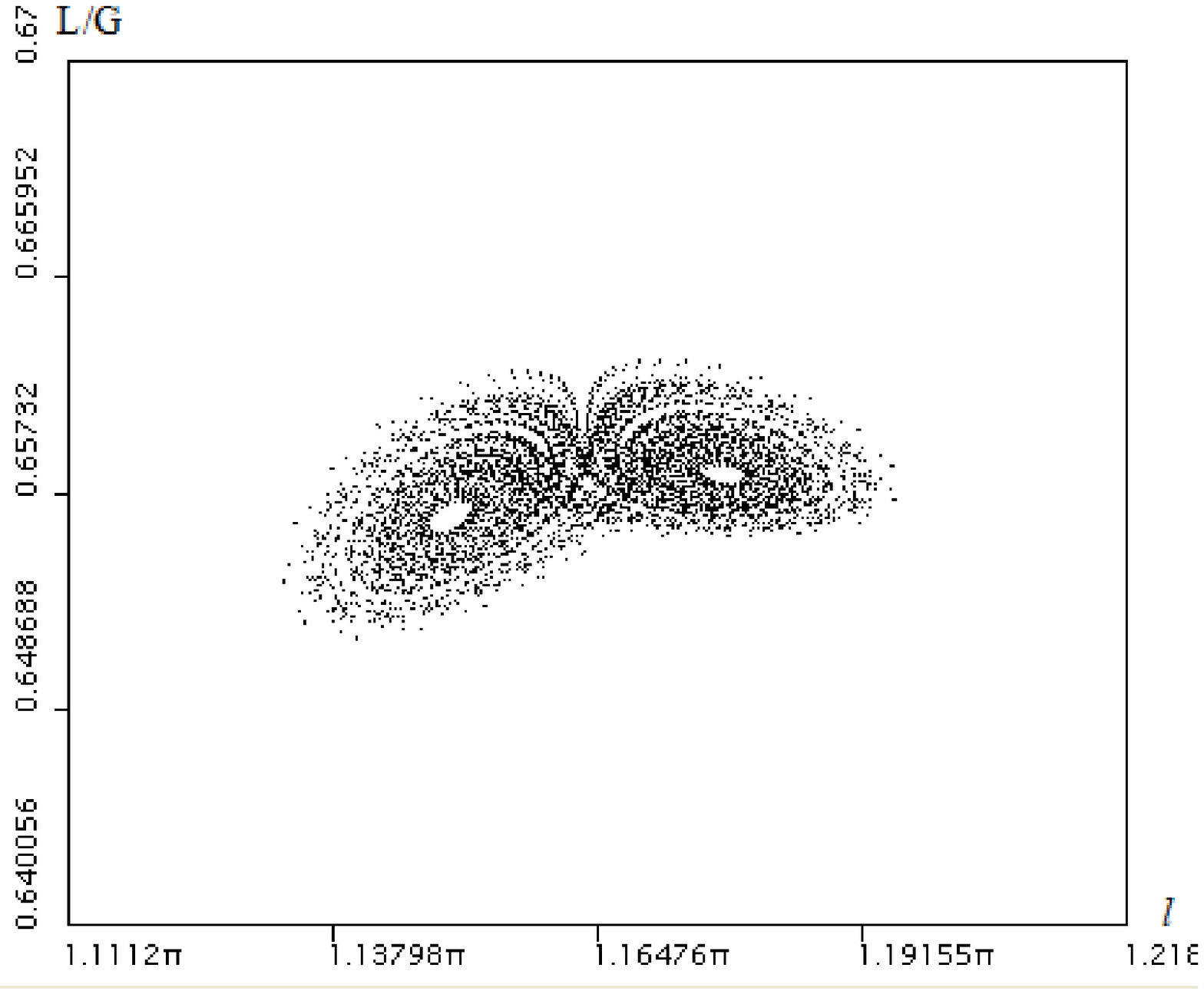,width=50mm}
 \\
d) E = 748.5 & e) E=750 & f) E=750 \\
\end{tabular}
\caption{{\footnotesize The main stages of evolution of the Lorenz-like attractor
in the map $T_E$. Figs. a) and f) show iterations of some starting point,
and Figs. b)--g) show unstable manifolds %(separatrices)
of the fixed point $O$.}}
\label{fig3-2a}
\end{figure}

4) Fig.~\ref{fig3-2a} illustrates the main stages of evolution of a discrete Lorenz attractor in the map ${\cal F}_{0E}$
for the parameter $E$ growing from $E = 747$ to $E=E^* = 750$.

Initially the attractor is a stable fixed point $O$, Fig.~\ref{fig3-2a}~a). Then,
at $E = E_1 = 747.61$, it undergoes a period-doubling bifurcation, and the stable cycle $P=(p_1,p_2)$
of period two becomes an attractor, Fig.~\ref{fig3-2a}~b).
%(см. также рис.\ref{scen1-1}(a)).
At $E=E_2 =748.4395$ the ``homoclinic figure-eight-butterfly'' of the unstable
manifolds (separatrices)
of the saddle $O$ is created, Fig.~\ref{fig3-2a}~c), which then gives rise
to a saddle-type closed invariant curve $L = (L_1,L_2)$ of period two (where ${\cal F}_{0E}(L_1)=L_2, {\cal F}_{0E}(L_2)=L_1$),
the curves $L_1$ and $L_2$ surround the point $p_1$ and $p_2$, respectively.
At the same time, the unstable separatrices of $O$ are rebuilt and now, for
$E_2<E<E_3$, the left (right) one is coiled around the right (left) point of the cycle $P$, Fig.~\ref{fig3-2a}~d).
Moreover, together with the closed period-$2$ invariant curve $L$, an invariant limit
set $\Omega$ is born here, \cite{Sh80}, which is not attracting yet.
As the numerical calculations show, for $E=E_3\sim 748.97$ the separatrices ``lie'' on the stable
manifold of the curve $L$ and then leave it. Almost immediately after that, at $E = E_4 \sim 748.98$, the
period-$2$ cycle $P$ sharply loses stability
under a subcritical torus-birth bifurcation: the closed invariant curve $L$ merges with the cycle $P$, after that the cycle becomes a saddle and the curve disappears.
The value of $E = E_4$ is the bifurcation moment of the creation of strange attractor -- the invariant set $\Omega$
becomes attracting. Even for the parameter values close to $E=E_3$ (and $E>E_3$) the
separatrices start to unwind,
see Fig.~\ref{fig3-2a}~e) and their configuration becomes similar to the Lorenzian one, which also
applies to the phase portrait, see Fig.~\ref{fig3-2a}~f). \\

\begin{figure}[htb]
\begin{tabular}{cc}
\psfig{file=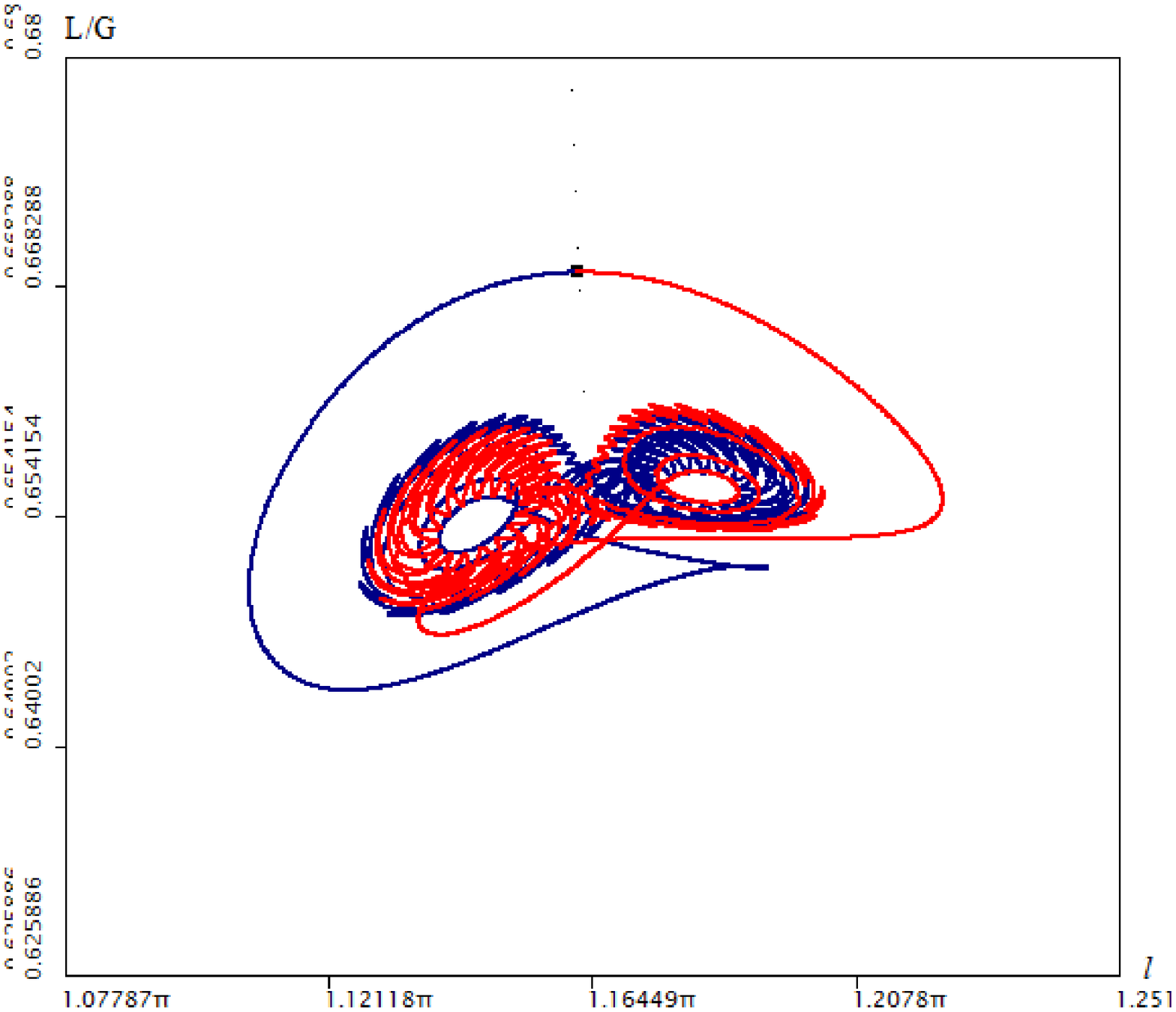,width=75mm} &
\psfig{file=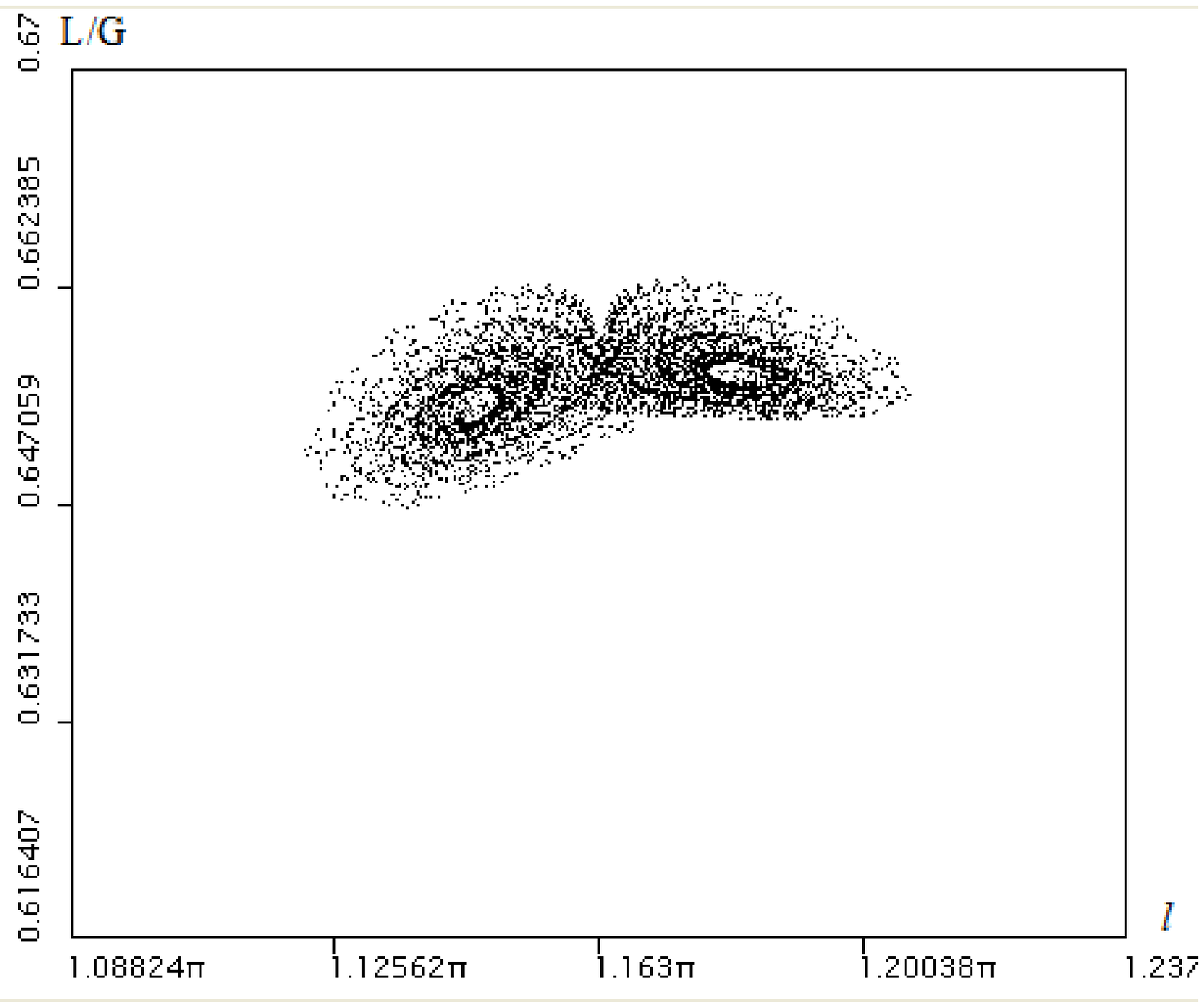,width=75mm} \\
a)  & b)
\end{tabular}
\caption{}
\label{fig3-2b}
\end{figure}

Fig.~\ref{fig3-2b} shows the behavior of (a) manifolds $W^u(O^*)$ and  b) iterations of the points on the attractor $A^*$
of map ${\cal F}_{0E^*}$ (here $E=E^*=752$). This attractor is studied
in items 1)--3) above.

\begin{figure}[htb]
\begin{tabular}{ccc}
\psfig{file=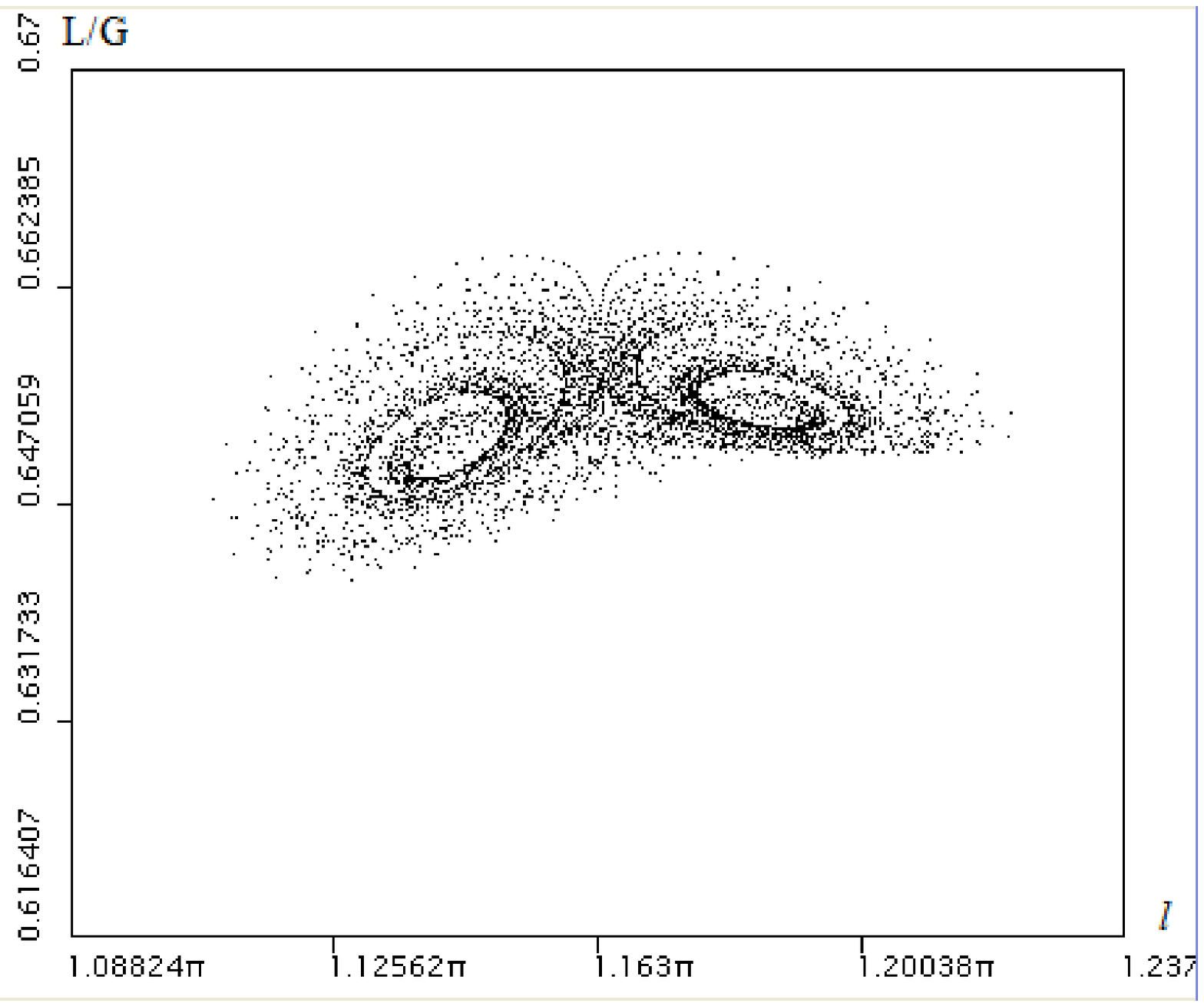,width=50mm}&
\psfig{file=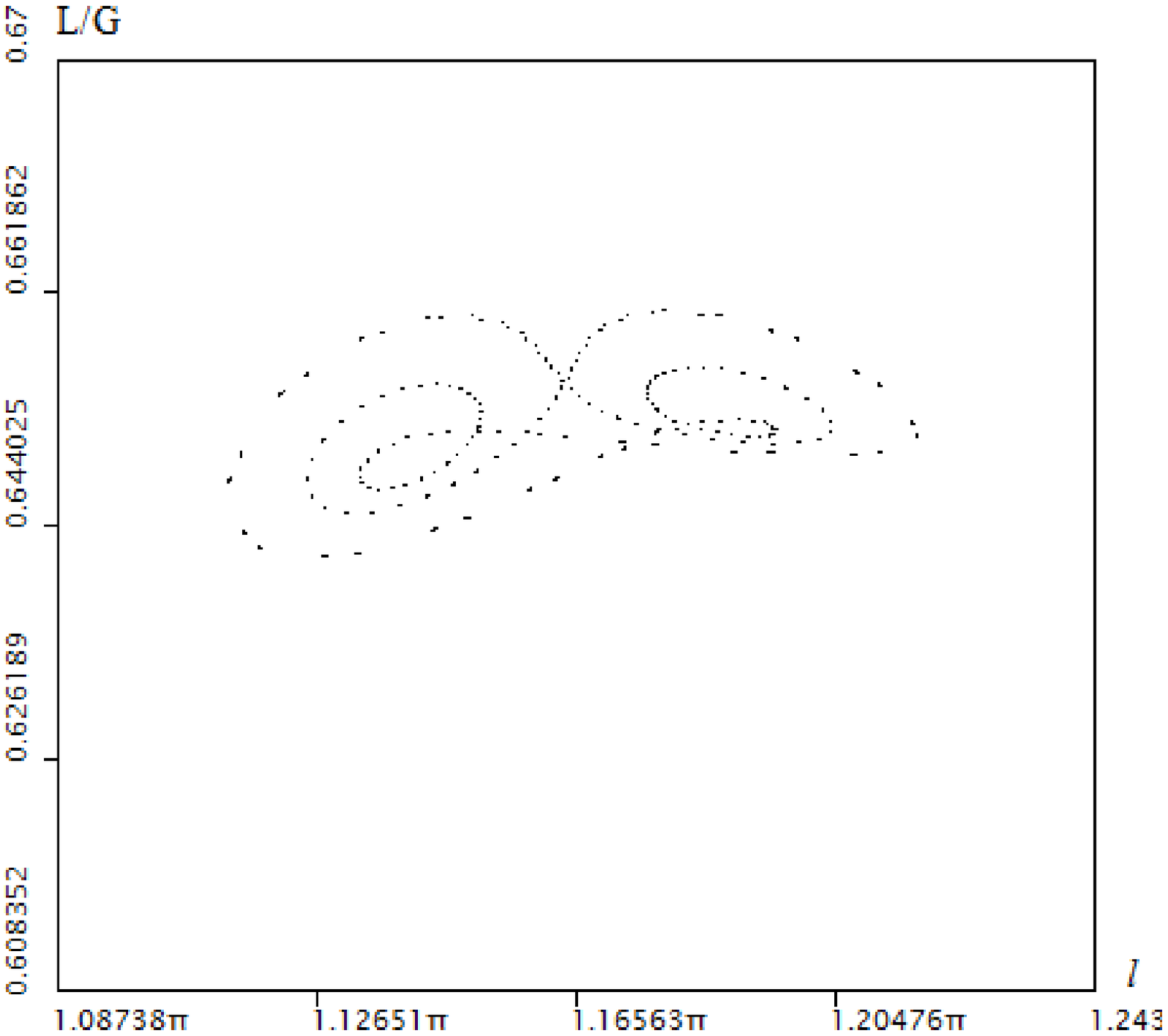, width=50mm} &
\psfig{file=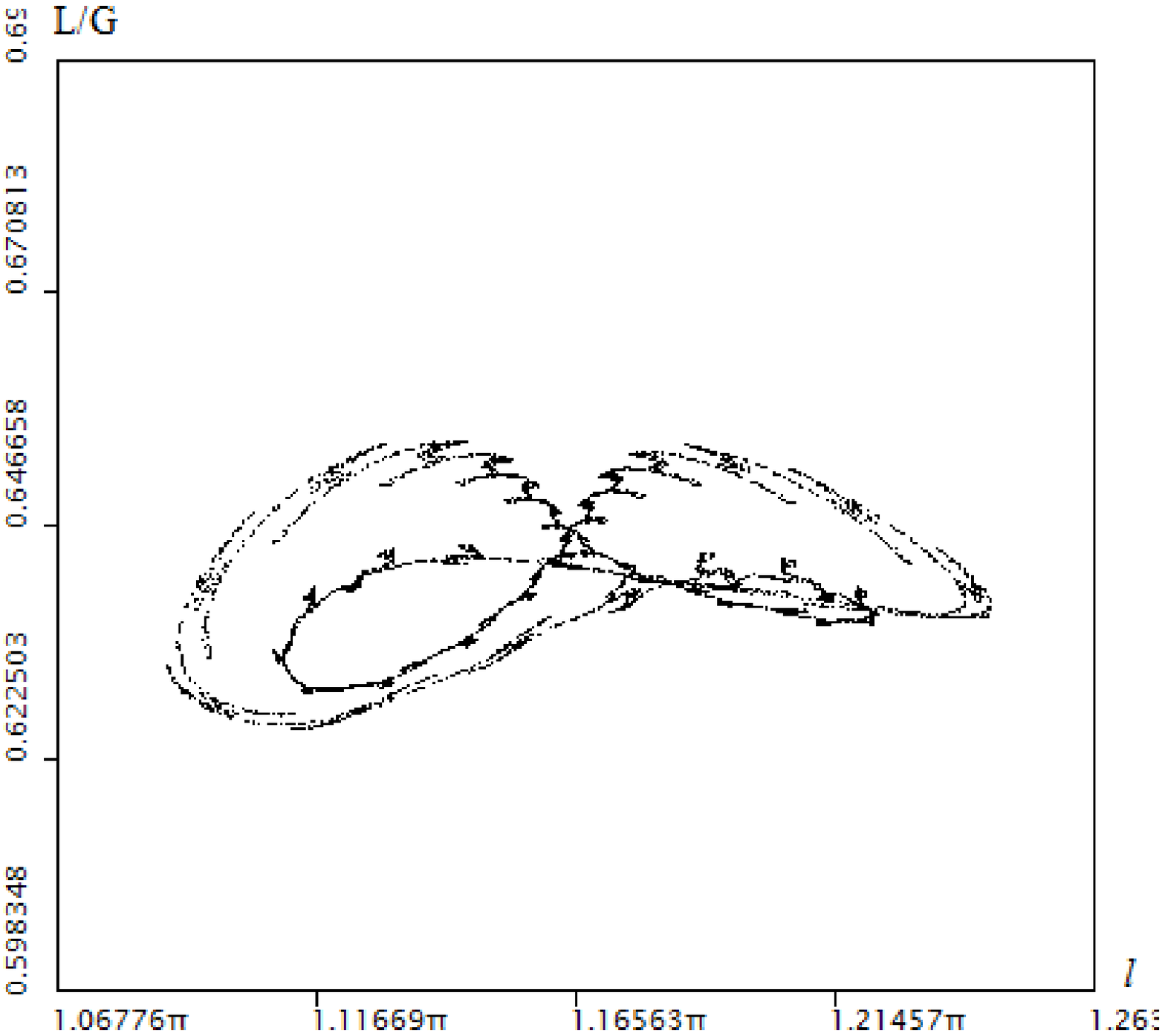,width=50mm} \\
a) E=754 & b) E=755 & c) E=765 \\
\psfig{file=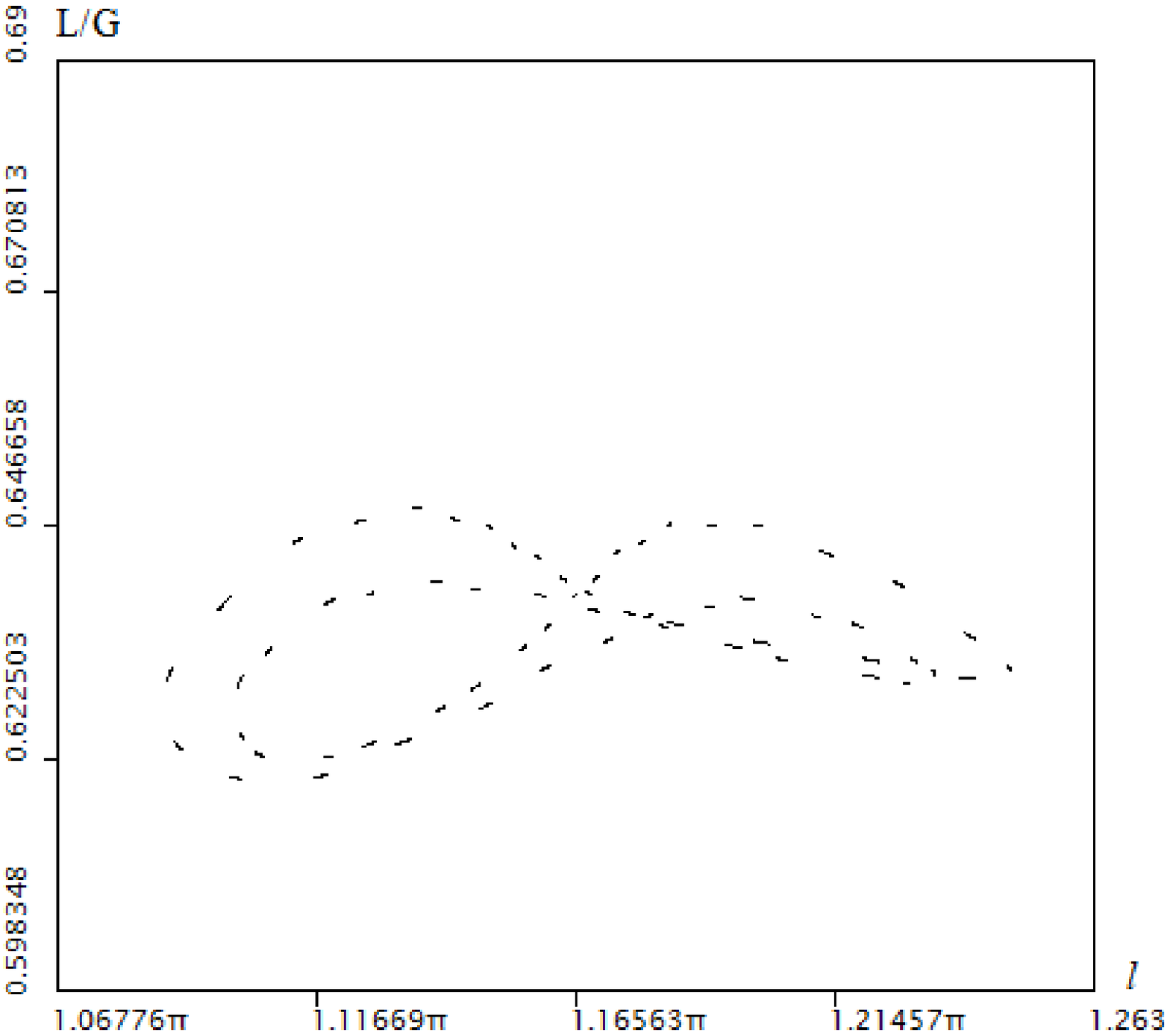,width=50mm} &
\psfig{file=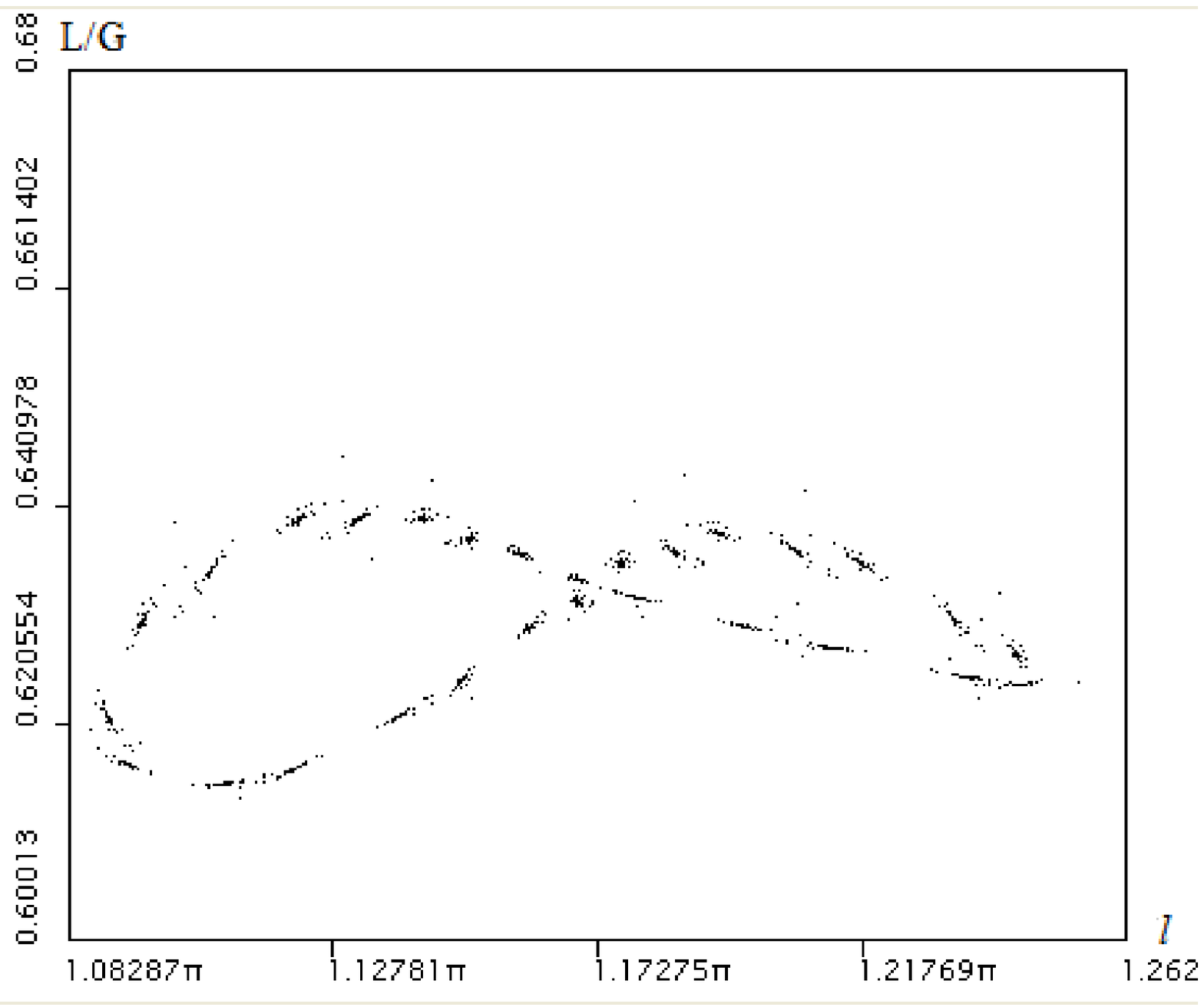, width=50mm} &
%
%d) E= 765 e) E=770 & f) E=775 \\
\psfig{file=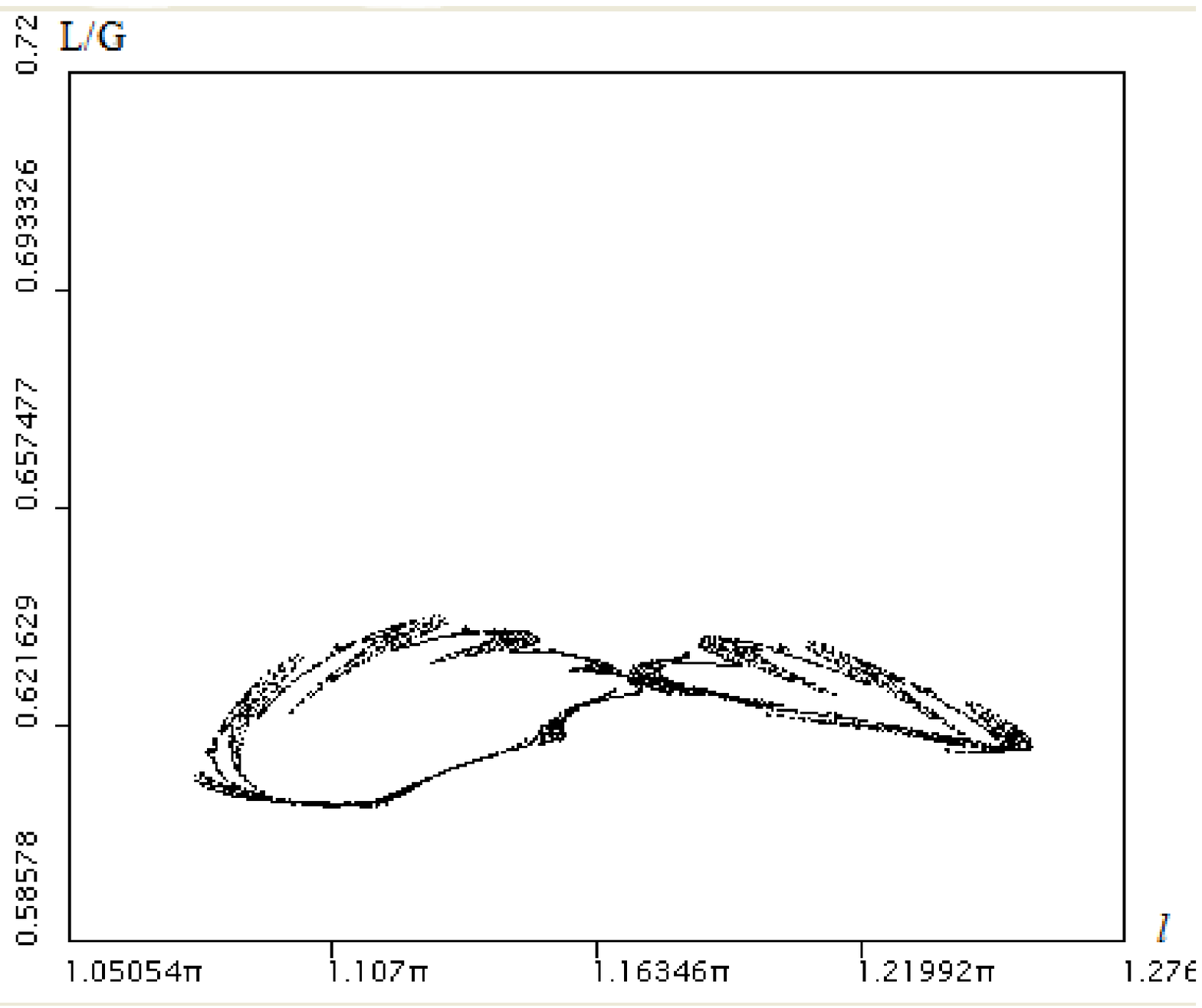,width=50mm} \\
d) E= 770 &e) E=775 & f) E=780 \\
\end{tabular}
\caption{{\footnotesize Certain stages of destruction of the discrete Lorenz attractor: a) this is definitely a discrete Lorenz attractor; b), d) and e) various resonant invariant curves ;
c) this attractor is very probably not of Lorenz type; f) a torus-chaos.} }
\label{fig3-3}
\end{figure}

Fig.~\ref{fig3-3} shows some stages of destruction of the discrete Lorenz attractor, which is
related to the appearance of resonant stable invariant curves , (b), (d) and (e), and the chaotic
regimes (torus-chaos),   (c) and (f). Note that for $E > 790$ nothing remains of the discrete Lorenz
attractor and the orbits run away from its neighborhood, tending to a new stable regime --
the spiral attractor, observed in~\cite{GGK12,KJSS12}.  %\\~\\

\section{Conclusion.}

The Lorenz attractors for flows play a special role in the theory of dynamical chaos. Until recently, these and hyperbolic attractors were the only ones which were classified
as ``genuine'' strange
attractors, which, in particular, do not allow the appearance of stable periodic orbits under small
perturbations. After publication of the paper \cite{TS98} by Turaev and Shilnikov the situation changed drastically. They not only provided an example of a wild hyperbolic spiral
attractor that must be regarded as a genuine strange attractor but also introduced a new class of {\em pseudo-hyperbolic attractors}. Thus, a new trend related to the study of such strange attractors appeared in the theory of dynamical chaos.

Discrete Lorenz attractors should be considered as very interesting examples of such genuine attractors. However, unlike the flow Lorenz attractors, their mathematical theory has not been constructed yet. Although, the basic elements of this theory already exist, \cite{TS98,TS08,GOST05,GGOT13}, and ``actively work''.

In particular, the example of a discrete Lorenz attractor found in this paper for the model of Celtic stone is, as we know, the first one for models from applications. However, we are sure that discrete Lorenz attractors exist in other models, one needs to make a more detailed search. Especially, since such attractors are not quite exotic, they could appear in dynamical systems (e.g.
 in three-dimensional maps) as a result of simple and universal
 bifurcation scenarios \cite{GGS12,GGKT14}.

{\bf Aknowledgements.} The author are grateful to  M. Malkin and D.Turaev for fruitful discussions and useful comments.
This work has been partially supported by the Russian Scientific Foundation Grant 14-41-00044 and the RFBR grant 13-01-00589.
Section~3 is carried out by the RSciF-grant (project No.14-12-00811).\\~\\
%This work was supported by the RFBR grants  13-01-00589 and
%13-01-97028-povolzhye,.

\end{document}